\documentclass{article}
\long\def\ig#1{\relax}
\ig{Thanks to Roberto Minio for this def'n.  Compare the def'n of
\comment in AMSTeX.}

\newcount \coefa
\newcount \coefb
\newcount \coefc
\newcount\tempcounta
\newcount\tempcountb
\newcount\tempcountc
\newcount\tempcountd
\newcount\xext
\newcount\yext
\newcount\xoff
\newcount\yoff
\newcount\gap%
\newcount\arrowtypea
\newcount\arrowtypeb
\newcount\arrowtypec
\newcount\arrowtyped
\newcount\arrowtypee
\newcount\height
\newcount\width
\newcount\xpos
\newcount\ypos
\newcount\run
\newcount\rise
\newcount\arrowlength
\newcount\halflength
\newcount\arrowtype
\newdimen\tempdimen
\newdimen\xlen
\newdimen\ylen
\newsavebox{\tempboxa}%
\newsavebox{\tempboxb}%
\newsavebox{\tempboxc}%

\makeatletter
\def\settypes(#1,#2,#3){\arrowtypea#1 \arrowtypeb#2 \arrowtypec#3}
\def\settoheight#1#2{\setbox\@tempboxa\hbox{#2}#1\ht\@tempboxa\relax}%
\def\settodepth#1#2{\setbox\@tempboxa\hbox{#2}#1\dp\@tempboxa\relax}%
\def\settokens[#1`#2`#3`#4]{%
     \def\tokena{#1}\def\tokenb{#2}\def\tokenc{#3}\def\tokend{#4}}
\def\setsqparms[#1`#2`#3`#4;#5`#6]{%
\arrowtypea #1
\arrowtypeb #2
\arrowtypec #3
\arrowtyped #4
\width #5
\height #6
}
\def\setpos(#1,#2){\xpos=#1 \ypos#2}

\def\bfig{\begin{picture}(\xext,\yext)(\xoff,\yoff)}
\def\efig{\end{picture}}

\def\putbox(#1,#2)#3{\put(#1,#2){\makebox(0,0){$#3$}}}

\def\settriparms[#1`#2`#3;#4]{\settripairparms[#1`#2`#3`1`1;#4]}%

\def\settripairparms[#1`#2`#3`#4`#5;#6]{%
\arrowtypea #1
\arrowtypeb #2
\arrowtypec #3
\arrowtyped #4
\arrowtypee #5
\width #6
\height #6
}

\def\resetparms{\settripairparms[1`1`1`1`1;500]\width 500}%default values%

\resetparms

\def\mvector(#1,#2)#3{%%
\put(0,0){\vector(#1,#2){#3}}%
\put(0,0){\vector(#1,#2){30}}%
}
\def\evector(#1,#2)#3{{%%
\arrowlength #3
\put(0,0){\vector(#1,#2){\arrowlength}}%
\advance \arrowlength by-30
\put(0,0){\vector(#1,#2){\arrowlength}}%
}}

\def\horsize#1#2{%
\settowidth{\tempdimen}{$#2$}%
#1=\tempdimen
\divide #1 by\unitlength
}

\def\vertsize#1#2{%
\settoheight{\tempdimen}{$#2$}%
#1=\tempdimen
\settodepth{\tempdimen}{$#2$}%
\advance #1 by\tempdimen
\divide #1 by\unitlength
}

\def\vertadjust[#1`#2`#3]{%
\vertsize{\tempcounta}{#1}%
\vertsize{\tempcountb}{#2}%
\ifnum \tempcounta<\tempcountb \tempcounta=\tempcountb \fi
\divide\tempcounta by2
\vertsize{\tempcountb}{#3}%
\ifnum \tempcountb>0 \advance \tempcountb by20 \fi
\ifnum \tempcounta<\tempcountb \tempcounta=\tempcountb \fi
}

\def\horadjust[#1`#2`#3]{%
\horsize{\tempcounta}{#1}%
\horsize{\tempcountb}{#2}%
\ifnum \tempcounta<\tempcountb \tempcounta=\tempcountb \fi
\divide\tempcounta by20
\horsize{\tempcountb}{#3}%
\ifnum \tempcountb>0 \advance \tempcountb by60 \fi
\ifnum \tempcounta<\tempcountb \tempcounta=\tempcountb \fi
}

\ig{ In this procedure, #1 is the paramater that sticks out all the way,
#2 sticks out the least and #3 is a label sticking out half way.  #4 is
the amount of the offset.}

\def\sladjust[#1`#2`#3]#4{%
\tempcountc=#4
\horsize{\tempcounta}{#1}%
\divide \tempcounta by2
\horsize{\tempcountb}{#2}%
\divide \tempcountb by2
\advance \tempcountb by-\tempcountc
\ifnum \tempcounta<\tempcountb \tempcounta=\tempcountb\fi
\divide \tempcountc by2
\horsize{\tempcountb}{#3}%
\advance \tempcountb by-\tempcountc
\ifnum \tempcountb>0 \advance \tempcountb by80\fi
\ifnum \tempcounta<\tempcountb \tempcounta=\tempcountb\fi
\advance\tempcounta by20
}

\def\putvector(#1,#2)(#3,#4)#5#6{{%
\xpos=#1
\ypos=#2
\run=#3
\rise=#4
\arrowlength=#5
\arrowtype=#6
\ifnum \arrowtype<0
    \ifnum \run=0
        \advance \ypos by-\arrowlength
    \else
        \tempcounta \arrowlength
        \multiply \tempcounta by\rise
        \divide \tempcounta by\run
        \ifnum\run>0
            \advance \xpos by\arrowlength
            \advance \ypos by\tempcounta
        \else
            \advance \xpos by-\arrowlength
            \advance \ypos by-\tempcounta
        \fi
    \fi
    \multiply \arrowtype by-1
    \multiply \rise by-1
    \multiply \run by-1
\fi
\ifnum \arrowtype=1
    \put(\xpos,\ypos){\vector(\run,\rise){\arrowlength}}%
\else\ifnum \arrowtype=2
    \put(\xpos,\ypos){\mvector(\run,\rise)\arrowlength}%
\else\ifnum\arrowtype=3
    \put(\xpos,\ypos){\evector(\run,\rise){\arrowlength}}%
\fi\fi\fi
}}

\def\putsplitvector(#1,#2)#3#4{%%
\xpos #1
\ypos #2
\arrowtype #4
\halflength #3
\arrowlength #3
\gap 140
\advance \halflength by-\gap
\divide \halflength by2
\ifnum \arrowtype=1
    \put(\xpos,\ypos){\line(0,-1){\halflength}}%
    \advance\ypos by-\halflength
    \advance\ypos by-\gap
    \put(\xpos,\ypos){\vector(0,-1){\halflength}}%
\else\ifnum \arrowtype=2
    \put(\xpos,\ypos){\line(0,-1)\halflength}%
    \put(\xpos,\ypos){\vector(0,-1)3}%
    \advance\ypos by-\halflength
    \advance\ypos by-\gap
    \put(\xpos,\ypos){\vector(0,-1){\halflength}}%
\else\ifnum\arrowtype=3
    \put(\xpos,\ypos){\line(0,-1)\halflength}%
    \advance\ypos by-\halflength
    \advance\ypos by-\gap
    \put(\xpos,\ypos){\evector(0,-1){\halflength}}%
\else\ifnum \arrowtype=-1
    \advance \ypos by-\arrowlength
    \put(\xpos,\ypos){\line(0,1){\halflength}}%
    \advance\ypos by\halflength
    \advance\ypos by\gap
    \put(\xpos,\ypos){\vector(0,1){\halflength}}%
\else\ifnum \arrowtype=-2
    \advance \ypos by-\arrowlength
    \put(\xpos,\ypos){\line(0,1)\halflength}%
    \put(\xpos,\ypos){\vector(0,1)3}%
    \advance\ypos by\halflength
    \advance\ypos by\gap
    \put(\xpos,\ypos){\vector(0,1){\halflength}}%
\else\ifnum\arrowtype=-3
    \advance \ypos by-\arrowlength
    \put(\xpos,\ypos){\line(0,1)\halflength}%
    \advance\ypos by\halflength
    \advance\ypos by\gap
    \put(\xpos,\ypos){\evector(0,1){\halflength}}%
\fi\fi\fi\fi\fi\fi
}

\def\putmorphism(#1)(#2,#3)[#4`#5`#6]#7#8#9{{%
\run #2
\rise #3
\ifnum\rise=0
  \puthmorphism(#1)[#4`#5`#6]{#7}{#8}{#9}%
\else\ifnum\run=0
  \putvmorphism(#1)[#4`#5`#6]{#7}{#8}{#9}%
\else
\setpos(#1)%
\arrowlength #7
\arrowtype #8
\ifnum\run=0
\else\ifnum\rise=0
\else
\ifnum\run>0
    \coefa=1
\else
   \coefa=-1
\fi
\ifnum\arrowtype>0
   \coefb=0
   \coefc=-1
\else
   \coefb=\coefa
   \coefc=1
   \arrowtype=-\arrowtype
\fi
\width=2
\multiply \width by\run
\divide \width by\rise
\ifnum \width<0  \width=-\width\fi
\advance\width by60
\if l#9 \width=-\width\fi
\putbox(\xpos,\ypos){#4}%            %node 1
{\multiply \coefa by\arrowlength%      %node 2
\advance\xpos by\coefa
\multiply \coefa by\rise
\divide \coefa by\run
\advance \ypos by\coefa
\putbox(\xpos,\ypos){#5} }%
{\multiply \coefa by\arrowlength%      %label
\divide \coefa by2
\advance \xpos by\coefa
\advance \xpos by\width
\multiply \coefa by\rise
\divide \coefa by\run
\advance \ypos by\coefa
\if l#9%
   \put(\xpos,\ypos){\makebox(0,0)[r]{$#6$}}%
\else\if r#9%
   \put(\xpos,\ypos){\makebox(0,0)[l]{$#6$}}%
\fi\fi }%
{\multiply \rise by-\coefc%             %arrow
\multiply \run by-\coefc
\multiply \coefb by\arrowlength
\advance \xpos by\coefb
\multiply \coefb by\rise
\divide \coefb by\run
\advance \ypos by\coefb
\multiply \coefc by70
\advance \ypos by\coefc
\multiply \coefc by\run
\divide \coefc by\rise
\advance \xpos by\coefc
\multiply \coefa by140
\multiply \coefa by\run
\divide \coefa by\rise
\advance \arrowlength by\coefa
\ifnum \arrowtype=1
   \put(\xpos,\ypos){\vector(\run,\rise){\arrowlength}}%
\else\ifnum\arrowtype=2
   \put(\xpos,\ypos){\mvector(\run,\rise){\arrowlength}}%
\else\ifnum\arrowtype=3
   \put(\xpos,\ypos){\evector(\run,\rise){\arrowlength}}%
\fi\fi\fi}\fi\fi\fi\fi}}

\def\puthmorphism(#1,#2)[#3`#4`#5]#6#7#8{{%
\xpos #1
\ypos #2
\width #6
\arrowlength #6
\putbox(\xpos,\ypos){#3\vphantom{#4}}%
{\advance \xpos by\arrowlength
\putbox(\xpos,\ypos){\vphantom{#3}#4}}%
\horsize{\tempcounta}{#3}%
\horsize{\tempcountb}{#4}%
\divide \tempcounta by2
\divide \tempcountb by2
\advance \tempcounta by30
\advance \tempcountb by30
\advance \xpos by\tempcounta
\advance \arrowlength by-\tempcounta
\advance \arrowlength by-\tempcountb
\putvector(\xpos,\ypos)(1,0){\arrowlength}{#7}%
\divide \arrowlength by2
\advance \xpos by\arrowlength
\vertsize{\tempcounta}{#5}%
\divide\tempcounta by2
\advance \tempcounta by20
\if a#8 %
   \advance \ypos by\tempcounta
   \putbox(\xpos,\ypos){#5}%
\else
   \advance \ypos by-\tempcounta
   \putbox(\xpos,\ypos){#5}%
\fi}}

\def\putvmorphism(#1,#2)[#3`#4`#5]#6#7#8{{%
\xpos #1
\ypos #2
\arrowlength #6
\arrowtype #7
\settowidth{\xlen}{$#5$}%
\putbox(\xpos,\ypos){#3}%
{\advance \ypos by-\arrowlength
\putbox(\xpos,\ypos){#4}}%
{\advance\arrowlength by-140
\advance \ypos by-70
\ifdim\xlen>0pt
   \if m#8%
      \putsplitvector(\xpos,\ypos){\arrowlength}{\arrowtype}%
   \else
      \putvector(\xpos,\ypos)(0,-1){\arrowlength}{\arrowtype}%
   \fi
\else
   \putvector(\xpos,\ypos)(0,-1){\arrowlength}{\arrowtype}%
\fi}%
\ifdim\xlen>0pt
   \divide \arrowlength by2
   \advance\ypos by-\arrowlength
   \if l#8%
      \advance \xpos by-40
      \put(\xpos,\ypos){\makebox(0,0)[r]{$#5$}}%
   \else\if r#8%
      \advance \xpos by40
      \put(\xpos,\ypos){\makebox(0,0)[l]{$#5$}}%
   \else
      \putbox(\xpos,\ypos){#5}%
   \fi\fi
\fi
}}

\def\topadjust[#1`#2`#3]{%
\yoff=10
\vertadjust[#1`#2`{#3}]%
\advance \yext by\tempcounta
\advance \yext by 10
}
\def\botadjust[#1`#2`#3]{%
\vertadjust[#1`#2`{#3}]%
\advance \yext by\tempcounta
\advance \yoff by-\tempcounta
}
\def\leftadjust[#1`#2`#3]{%
\xoff=0
\horadjust[#1`#2`{#3}]%
\advance \xext by\tempcounta
\advance \xoff by-\tempcounta
}
\def\rightadjust[#1`#2`#3]{%
\horadjust[#1`#2`{#3}]%
\advance \xext by\tempcounta
}
\def\rightsladjust[#1`#2`#3]{%
\sladjust[#1`#2`{#3}]{\width}%
\advance \xext by\tempcounta
}
\def\leftsladjust[#1`#2`#3]{%
\xoff=0
\sladjust[#1`#2`{#3}]{\width}%
\advance \xext by\tempcounta
\advance \xoff by-\tempcounta
}
\def\adjust[#1`#2;#3`#4;#5`#6;#7`#8]{%
\topadjust[#1``{#2}]
\leftadjust[#3``{#4}]
\rightadjust[#5``{#6}]
\botadjust[#7``{#8}]}

\def\putsquarep<#1>(#2)[#3;#4`#5`#6`#7]{{%
\setsqparms[#1]%
\setpos(#2)%
\settokens[#3]%
\puthmorphism(\xpos,\ypos)[\tokenc`\tokend`{#7}]{\width}{\arrowtyped}b%
\advance\ypos by \height
\puthmorphism(\xpos,\ypos)[\tokena`\tokenb`{#4}]{\width}{\arrowtypea}a%
\putvmorphism(\xpos,\ypos)[``{#5}]{\height}{\arrowtypeb}l%
\advance\xpos by \width
\putvmorphism(\xpos,\ypos)[``{#6}]{\height}{\arrowtypec}r%
}}

\def\putsquare{\@ifnextchar <{\putsquarep}{\putsquarep%
   <\arrowtypea`\arrowtypeb`\arrowtypec`\arrowtyped;\width`\height>}}
\def\square{\@ifnextchar< {\squarep}{\squarep
   <\arrowtypea`\arrowtypeb`\arrowtypec`\arrowtyped;\width`\height>}}
\def\squarep<#1>[#2`#3`#4`#5;#6`#7`#8`#9]{{%          %     #2------>#3
\setsqparms[#1]%                                      %      |       |
\xext=\width                                          %      |       |
\yext=\height                                         %    #7|       |#8
\topadjust[#2`#3`{#6}]%                               %      |       |
\botadjust[#4`#5`{#9}]%                               %      |       |
\leftadjust[#2`#4`{#7}]%                              %
\rightadjust[#3`#5`{#8}]%                             %     #4------>#5
\begin{picture}(\xext,\yext)(\xoff,\yoff)%                      #9
\putsquarep<\arrowtypea`\arrowtypeb`\arrowtypec`\arrowtyped;\width`\height>%
(0,0)[#2`#3`#4`#5;#6`#7`#8`{#9}]%
\end{picture}%
}}

\def\putptrianglep<#1>(#2,#3)[#4`#5`#6;#7`#8`#9]{{%
\settriparms[#1]%
\xpos=#2 \ypos=#3
\advance\ypos by \height
\puthmorphism(\xpos,\ypos)[#4`#5`{#7}]{\height}{\arrowtypea}a%
\putvmorphism(\xpos,\ypos)[`#6`{#8}]{\height}{\arrowtypeb}l%
\advance\xpos by\height
\putmorphism(\xpos,\ypos)(-1,-1)[``{#9}]{\height}{\arrowtypec}r%
}}

\def\putptriangle{\@ifnextchar <{\putptrianglep}{\putptrianglep
   <\arrowtypea`\arrowtypeb`\arrowtypec;\height>}}
\def\ptriangle{\@ifnextchar <{\ptrianglep}{\ptrianglep
   <\arrowtypea`\arrowtypeb`\arrowtypec;\height>}}

\def\ptrianglep<#1>[#2`#3`#4;#5`#6`#7]{{%%       #5
\settriparms[#1]%
\width=\height                         %      #2----->#3
\xext=\width                           %      |      /
\yext=\width                           %      |     /
\topadjust[#2`#3`{#5}]%                %    #6|    /#7
\botadjust[#3``]%                      %      |   /
\leftadjust[#2`#4`{#6}]%               %      |  /
\rightsladjust[#3`#4`{#7}]%            %
\begin{picture}(\xext,\yext)(\xoff,\yoff)%    #4
\putptrianglep<\arrowtypea`\arrowtypeb`\arrowtypec;\height>%
(0,0)[#2`#3`#4;#5`#6`{#7}]%
\end{picture}%
}}

\def\putqtrianglep<#1>(#2,#3)[#4`#5`#6;#7`#8`#9]{{%
\settriparms[#1]%
\xpos=#2 \ypos=#3
\advance\ypos by\height
\puthmorphism(\xpos,\ypos)[#4`#5`{#7}]{\height}{\arrowtypea}a%
\putmorphism(\xpos,\ypos)(1,-1)[``{#8}]{\height}{\arrowtypeb}l%
\advance\xpos by\height
\putvmorphism(\xpos,\ypos)[`#6`{#9}]{\height}{\arrowtypec}r%
}}

\def\putqtriangle{\@ifnextchar <{\putqtrianglep}{\putqtrianglep
   <\arrowtypea`\arrowtypeb`\arrowtypec;\height>}}
\def\qtriangle{\@ifnextchar <{\qtrianglep}{\qtrianglep
   <\arrowtypea`\arrowtypeb`\arrowtypec;\height>}}

\def\qtrianglep<#1>[#2`#3`#4;#5`#6`#7]{{%%
\settriparms[#1]%                                  #5
\width=\height                         %        #2----->#3
\xext=\width                           %         \      |
\yext=\height                          %          \     |
\topadjust[#2`#3`{#5}]%                %         #6\    |#7
\botadjust[#4``]%                      %            \   |
\leftsladjust[#2`#4`{#6}]%             %             \  |
\rightadjust[#3`#4`{#7}]%              %
\begin{picture}(\xext,\yext)(\xoff,\yoff)%             #4
\putqtrianglep<\arrowtypea`\arrowtypeb`\arrowtypec;\height>%
(0,0)[#2`#3`#4;#5`#6`{#7}]%
\end{picture}%
}}

\def\putdtrianglep<#1>(#2,#3)[#4`#5`#6;#7`#8`#9]{{%
\settriparms[#1]%
\xpos=#2 \ypos=#3
\puthmorphism(\xpos,\ypos)[#5`#6`{#9}]{\height}{\arrowtypec}b%
\advance\xpos by \height \advance\ypos by\height
\putmorphism(\xpos,\ypos)(-1,-1)[``{#7}]{\height}{\arrowtypea}l%
\putvmorphism(\xpos,\ypos)[#4``{#8}]{\height}{\arrowtypeb}r%
}}

\def\putdtriangle{\@ifnextchar <{\putdtrianglep}{\putdtrianglep
   <\arrowtypea`\arrowtypeb`\arrowtypec;\height>}}
\def\dtriangle{\@ifnextchar <{\dtrianglep}{\dtrianglep
   <\arrowtypea`\arrowtypeb`\arrowtypec;\height>}}

\def\dtrianglep<#1>[#2`#3`#4;#5`#6`#7]{{%%
\settriparms[#1]%                                          #2
\width=\height                         %                  / |
\xext=\width                           %                 /  |
\yext=\height                          %              #5/   |#6
\topadjust[#2``]%                      %               /    |
\botadjust[#3`#4`{#7}]%                %              /     |
\leftsladjust[#3`#2`{#5}]%             %
\rightadjust[#2`#4`{#6}]%              %            #3----->#4
\begin{picture}(\xext,\yext)(\xoff,\yoff)%              #7
\putdtrianglep<\arrowtypea`\arrowtypeb`\arrowtypec;\height>%
(0,0)[#2`#3`#4;#5`#6`{#7}]%
\end{picture}%
}}

\def\putbtrianglep<#1>(#2,#3)[#4`#5`#6;#7`#8`#9]{{%
\settriparms[#1]%
\xpos=#2 \ypos=#3
\puthmorphism(\xpos,\ypos)[#5`#6`{#9}]{\height}{\arrowtypec}b%
\advance\ypos by\height
\putmorphism(\xpos,\ypos)(1,-1)[``{#8}]{\height}{\arrowtypeb}r%
\putvmorphism(\xpos,\ypos)[#4``{#7}]{\height}{\arrowtypea}l%
}}

\def\putbtriangle{\@ifnextchar <{\putbtrianglep}{\putbtrianglep
   <\arrowtypea`\arrowtypeb`\arrowtypec;\height>}}
\def\btriangle{\@ifnextchar <{\btrianglep}{\btrianglep
   <\arrowtypea`\arrowtypeb`\arrowtypec;\height>}}

\def\btrianglep<#1>[#2`#3`#4;#5`#6`#7]{{%%
\settriparms[#1]%                                     #2
\width=\height                         %              | \
\xext=\width                           %              |  \
\yext=\height                          %            #5|   \#6
\topadjust[#2``]%                      %              |    \
\botadjust[#3`#4`{#7}]%                %              |     \
\leftadjust[#2`#3`{#5}]%               %
\rightsladjust[#4`#2`{#6}]%            %              #3----->#4
\begin{picture}(\xext,\yext)(\xoff,\yoff)%                #7
\putbtrianglep<\arrowtypea`\arrowtypeb`\arrowtypec;\height>%
(0,0)[#2`#3`#4;#5`#6`{#7}]%
\end{picture}%
}}

\def\putAtrianglep<#1>(#2,#3)[#4`#5`#6;#7`#8`#9]{{%
\settriparms[#1]%
\xpos=#2 \ypos=#3
{\multiply \height by2
\puthmorphism(\xpos,\ypos)[#5`#6`{#9}]{\height}{\arrowtypec}b}%
\advance\xpos by\height \advance\ypos by\height
\putmorphism(\xpos,\ypos)(-1,-1)[#4``{#7}]{\height}{\arrowtypea}l%
\putmorphism(\xpos,\ypos)(1,-1)[``{#8}]{\height}{\arrowtypeb}r%
}}

\def\putAtriangle{\@ifnextchar <{\putAtrianglep}{\putAtrianglep
   <\arrowtypea`\arrowtypeb`\arrowtypec;\height>}}
\def\Atriangle{\@ifnextchar <{\Atrianglep}{\Atrianglep
   <\arrowtypea`\arrowtypeb`\arrowtypec;\height>}}

\def\Atrianglep<#1>[#2`#3`#4;#5`#6`#7]{{%%
\settriparms[#1]%                                 #2
\width=\height                         %         /   \
\xext=\width                           %        /     \
\yext=\height                          %     #5/       \#6
\topadjust[#2``]%                      %      /         \
\botadjust[#3`#4`{#7}]%                %     /           \
\multiply \xext by2 %                  %
\leftsladjust[#3`#2`{#5}]%             %   #3------------>#4
\rightsladjust[#4`#2`{#6}]%            %          #7
\begin{picture}(\xext,\yext)(\xoff,\yoff)%
\putAtrianglep<\arrowtypea`\arrowtypeb`\arrowtypec;\height>%
(0,0)[#2`#3`#4;#5`#6`{#7}]%
\end{picture}%
}}

\def\putAtrianglepairp<#1>(#2)[#3;#4`#5`#6`#7`#8]{{
\settripairparms[#1]%
\setpos(#2)%
\settokens[#3]%
\puthmorphism(\xpos,\ypos)[\tokenb`\tokenc`{#7}]{\height}{\arrowtyped}b%
\advance\xpos by\height
\advance\ypos by\height
\putmorphism(\xpos,\ypos)(-1,-1)[\tokena``{#4}]{\height}{\arrowtypea}l%
\putvmorphism(\xpos,\ypos)[``{#5}]{\height}{\arrowtypeb}m%
\putmorphism(\xpos,\ypos)(1,-1)[``{#6}]{\height}{\arrowtypec}r%
}}

\def\putAtrianglepair{\@ifnextchar <{\putAtrianglepairp}{\putAtrianglepairp%
   <\arrowtypea`\arrowtypeb`\arrowtypec`\arrowtyped`\arrowtypee;\height>}}
\def\Atrianglepair{\@ifnextchar <{\Atrianglepairp}{\Atrianglepairp%
   <\arrowtypea`\arrowtypeb`\arrowtypec`\arrowtyped`\arrowtypee;\height>}}

\def\Atrianglepairp<#1>[#2;#3`#4`#5`#6`#7]{{%
\settripairparms[#1]%
\settokens[#2]%
\width=\height
\xext=\width
\yext=\height
\topadjust[\tokena``]%
\vertadjust[\tokenb`\tokenc`{#6}]%                      %  #2a
\tempcountd=\tempcounta                       %           / | \
\vertadjust[\tokenc`\tokend`{#7}]%            %          /  |  \
\ifnum\tempcounta<\tempcountd                 %       #3/  #4   \#5
\tempcounta=\tempcountd\fi                    %        /    |    \
\advance \yext by\tempcounta                  %       /     |     \
\advance \yoff by-\tempcounta                 %
\multiply \xext by2 %                         %     #2b---->#2c---->#2d
\leftsladjust[\tokenb`\tokena`{#3}]%          %         #6     #7
\rightsladjust[\tokend`\tokena`{#5}]%
\begin{picture}(\xext,\yext)(\xoff,\yoff)%
\putAtrianglepairp
<\arrowtypea`\arrowtypeb`\arrowtypec`\arrowtyped`\arrowtypee;\height>%
(0,0)[#2;#3`#4`#5`#6`{#7}]%
\end{picture}%
}}

\def\putVtrianglep<#1>(#2,#3)[#4`#5`#6;#7`#8`#9]{{%
\settriparms[#1]%
\xpos=#2 \ypos=#3
\advance\ypos by\height
{\multiply\height by2
\puthmorphism(\xpos,\ypos)[#4`#5`{#7}]{\height}{\arrowtypea}a}%
\putmorphism(\xpos,\ypos)(1,-1)[`#6`{#8}]{\height}{\arrowtypeb}l%
\advance\xpos by\height
\advance\xpos by\height
\putmorphism(\xpos,\ypos)(-1,-1)[``{#9}]{\height}{\arrowtypec}r%
}}

\def\putVtriangle{\@ifnextchar <{\putVtrianglep}{\putVtrianglep
   <\arrowtypea`\arrowtypeb`\arrowtypec;\height>}}
\def\Vtriangle{\@ifnextchar <{\Vtrianglep}{\Vtrianglep
   <\arrowtypea`\arrowtypeb`\arrowtypec;\height>}}

\def\Vtrianglep<#1>[#2`#3`#4;#5`#6`#7]{{%%
\settriparms[#1]%                                      #5
\width=\height                         %        #2------------->#3
\xext=\width                           %         \             /
\yext=\height                          %          \           /
\topadjust[#2`#3`{#5}]%                %         #6\         /#7
\botadjust[#4``]%                      %            \       /
\multiply \xext by2 %                  %             \     /
\leftsladjust[#2`#3`{#6}]%             %
\rightsladjust[#3`#4`{#7}]%            %               #4
\begin{picture}(\xext,\yext)(\xoff,\yoff)%
\putVtrianglep<\arrowtypea`\arrowtypeb`\arrowtypec;\height>%
(0,0)[#2`#3`#4;#5`#6`{#7}]%
\end{picture}%
}}

\def\putVtrianglepairp<#1>(#2)[#3;#4`#5`#6`#7`#8]{{
\settripairparms[#1]%
\setpos(#2)%
\settokens[#3]%
\advance\ypos by\height
\putmorphism(\xpos,\ypos)(1,-1)[`\tokend`{#6}]{\height}{\arrowtypec}l%
\puthmorphism(\xpos,\ypos)[\tokena`\tokenb`{#4}]{\height}{\arrowtypea}a%
\advance\xpos by\height
\putvmorphism(\xpos,\ypos)[``{#7}]{\height}{\arrowtyped}m%
\advance\xpos by\height
\putmorphism(\xpos,\ypos)(-1,-1)[``{#8}]{\height}{\arrowtypee}r%
}}

\def\putVtrianglepair{\@ifnextchar <{\putVtrianglepairp}{\putVtrianglepairp%
    <\arrowtypea`\arrowtypeb`\arrowtypec`\arrowtyped`\arrowtypee;\height>}}
\def\Vtrianglepair{\@ifnextchar <{\Vtrianglepairp}{\Vtrianglepairp%
    <\arrowtypea`\arrowtypeb`\arrowtypec`\arrowtyped`\arrowtypee;\height>}}

\def\Vtrianglepairp<#1>[#2;#3`#4`#5`#6`#7]{{%
\settripairparms[#1]%
\settokens[#2]%                            #3      #4
\xext=\height                  %        #2a---->#2b---->#2c
\width=\height                 %         \      |      /
\yext=\height                  %          \     |     /
\vertadjust[\tokena`\tokenb`{#4}]%       #5\   #6    /#7
\tempcountd=\tempcounta        %            \   |   /
\vertadjust[\tokenb`\tokenc`{#5}]%           \  |  /
\ifnum\tempcounta<\tempcountd%
\tempcounta=\tempcountd\fi%                    #2d
\advance \yext by\tempcounta
\botadjust[\tokend``]%
\multiply \xext by2
\leftsladjust[\tokena`\tokend`{#6}]%
\rightsladjust[\tokenc`\tokend`{#7}]%
\begin{picture}(\xext,\yext)(\xoff,\yoff)%
\putVtrianglepairp
<\arrowtypea`\arrowtypeb`\arrowtypec`\arrowtyped`\arrowtypee;\height>%
(0,0)[#2;#3`#4`#5`#6`{#7}]%
\end{picture}%
}}

\def\putCtrianglep<#1>(#2,#3)[#4`#5`#6;#7`#8`#9]{{%
\settriparms[#1]%
\xpos=#2 \ypos=#3
\advance\ypos by\height
\putmorphism(\xpos,\ypos)(1,-1)[``{#9}]{\height}{\arrowtypec}l%
\advance\xpos by\height
\advance\ypos by\height
\putmorphism(\xpos,\ypos)(-1,-1)[#4`#5`{#7}]{\height}{\arrowtypea}l%
{\multiply\height by 2
\putvmorphism(\xpos,\ypos)[`#6`{#8}]{\height}{\arrowtypeb}r}%
}}

\def\putCtriangle{\@ifnextchar <{\putCtrianglep}{\putCtrianglep
    <\arrowtypea`\arrowtypeb`\arrowtypec;\height>}}
\def\Ctriangle{\@ifnextchar <{\Ctrianglep}{\Ctrianglep
    <\arrowtypea`\arrowtypeb`\arrowtypec;\height>}}

\def\Ctrianglep<#1>[#2`#3`#4;#5`#6`#7]{{%%
\settriparms[#1]%                                         #2
\width=\height                          %                / |
\xext=\width                            %               /  |
\yext=\height                           %            #5/   |
\multiply \yext by2 %                   %             /    |
\topadjust[#2``]%                       %            /     |
\botadjust[#4``]%                       %           v      |
\sladjust[#3`#2`{#5}]{\width}%          %          #3      |#6
\tempcountd=\tempcounta                 %           \      |
\sladjust[#3`#4`{#7}]{\width}%          %            \     |
\ifnum \tempcounta<\tempcountd          %           #7\    |
\tempcounta=\tempcountd\fi              %              \   |
\advance \xext by\tempcounta            %               \  |
\advance \xoff by-\tempcounta           %
\rightadjust[#2`#4`{#6}]%               %                 #4
\begin{picture}(\xext,\yext)(\xoff,\yoff)%
\putCtrianglep<\arrowtypea`\arrowtypeb`\arrowtypec;\height>%
(0,0)[#2`#3`#4;#5`#6`{#7}]%
\end{picture}%
}}

\def\putDtrianglep<#1>(#2,#3)[#4`#5`#6;#7`#8`#9]{{%
\settriparms[#1]%
\xpos=#2 \ypos=#3
\advance\xpos by\height \advance\ypos by\height
\putmorphism(\xpos,\ypos)(-1,-1)[``{#9}]{\height}{\arrowtypec}r%
\advance\xpos by-\height \advance\ypos by\height
\putmorphism(\xpos,\ypos)(1,-1)[`#5`{#8}]{\height}{\arrowtypeb}r%
{\multiply\height by 2
\putvmorphism(\xpos,\ypos)[#4`#6`{#7}]{\height}{\arrowtypea}l}%
}}

\def\putDtriangle{\@ifnextchar <{\putDtrianglep}{\putDtrianglep
    <\arrowtypea`\arrowtypeb`\arrowtypec;\height>}}
\def\Dtriangle{\@ifnextchar <{\Dtrianglep}{\Dtrianglep
   <\arrowtypea`\arrowtypeb`\arrowtypec;\height>}}

\def\Dtrianglep<#1>[#2`#3`#4;#5`#6`#7]{{%%
\settriparms[#1]%                                 #2
\width=\height                         %          | \
\xext=\height                          %          |  \
\yext=\height                          %          |   \#6
\multiply \yext by2 %                  %          |    \
\topadjust[#2``]%                      %          |     \
\botadjust[#4``]%                      %          |
\leftadjust[#2`#4`{#5}]%               %        #5|      #3
\sladjust[#3`#2`{#5}]{\height}%        %          |      /
\tempcountd=\tempcountd                %          |     /
\sladjust[#3`#4`{#7}]{\height}%        %          |    /#7
\ifnum \tempcounta<\tempcountd         %          |   /
\tempcounta=\tempcountd\fi             %          |  /
\advance \xext by\tempcounta           %
\begin{picture}(\xext,\yext)(\xoff,\yoff)%        #4
\putDtrianglep<\arrowtypea`\arrowtypeb`\arrowtypec;\height>%
(0,0)[#2`#3`#4;#5`#6`{#7}]%
\end{picture}%
}}

\def\setrecparms[#1`#2]{\width=#1 \height=#2}%
\def\recursep<#1`#2>[#3;#4`#5`#6`#7`#8]{{%
\width=#1 \height=#2
\settokens[#3]
\settowidth{\tempdimen}{$\tokena$}
\ifdim\tempdimen=0pt
  \savebox{\tempboxa}{\hbox{$\tokenb$}}%
  \savebox{\tempboxb}{\hbox{$\tokend$}}%
  \savebox{\tempboxc}{\hbox{$#6$}}%
\else
  \savebox{\tempboxa}{\hbox{$\hbox{$\tokena$}\times\hbox{$\tokenb$}$}}%
  \savebox{\tempboxb}{\hbox{$\hbox{$\tokena$}\times\hbox{$\tokend$}$}}%
  \savebox{\tempboxc}{\hbox{$\hbox{$\tokena$}\times\hbox{$#6$}$}}%
\fi
\ypos=\height
\divide\ypos by 2
\xpos=\ypos
\advance\xpos by \width
\xext=\xpos \yext=\height
\topadjust[#3`\usebox{\tempboxa}`{#4}]%
\botadjust[#5`\usebox{\tempboxb}`{#8}]%
\sladjust[\tokenc`\tokenb`{#5}]{\ypos}%
\tempcountd=\tempcounta
\sladjust[\tokenc`\tokend`{#5}]{\ypos}%
\ifnum \tempcounta<\tempcountd
\tempcounta=\tempcountd\fi
\advance \xext by\tempcounta
\advance \xoff by-\tempcounta
\rightadjust[\usebox{\tempboxa}`\usebox{\tempboxb}`\usebox{\tempboxc}]%
\bfig
\putCtrianglep<-1`1`1;\ypos>(0,0)[`\tokenc`;#5`#6`{#7}]%
\puthmorphism(\ypos,0)[\tokend`\usebox{\tempboxb}`{#8}]{\width}{-1}b%
\puthmorphism(\ypos,\height)[\tokenb`\usebox{\tempboxa}`{#4}]{\width}{-1}a%
\advance\ypos by \width
\putvmorphism(\ypos,\height)[``\usebox{\tempboxc}]{\height}1r%
\efig
}}

\def\recurse{\@ifnextchar <{\recursep}{\recursep<\width`\height>}}

\def\puttwohmorphisms(#1,#2)[#3`#4;#5`#6]#7#8#9{{%
\puthmorphism(#1,#2)[#3`#4`]{#7}0a
\ypos=#2
\advance\ypos by 20
\puthmorphism(#1,\ypos)[\phantom{#3}`\phantom{#4}`#5]{#7}{#8}a
\advance\ypos by -40
\puthmorphism(#1,\ypos)[\phantom{#3}`\phantom{#4}`#6]{#7}{#9}b
}}

\def\puttwovmorphisms(#1,#2)[#3`#4;#5`#6]#7#8#9{{%
\putvmorphism(#1,#2)[#3`#4`]{#7}0a
\xpos=#1
\advance\xpos by -20
\putvmorphism(\xpos,#2)[\phantom{#3}`\phantom{#4}`#5]{#7}{#8}l
\advance\xpos by 40
\putvmorphism(\xpos,#2)[\phantom{#3}`\phantom{#4}`#6]{#7}{#9}r
}}

\def\puthcoequalizer(#1)[#2`#3`#4;#5`#6`#7]#8#9{{%
\setpos(#1)%
\puttwohmorphisms(\xpos,\ypos)[#2`#3;#5`#6]{#8}11%
\advance\xpos by #8
\puthmorphism(\xpos,\ypos)[\phantom{#3}`#4`#7]{#8}1{#9}
}}

\def\putvcoequalizer(#1)[#2`#3`#4;#5`#6`#7]#8#9{{%
\setpos(#1)%
\puttwovmorphisms(\xpos,\ypos)[#2`#3;#5`#6]{#8}11%
\advance\ypos by -#8
\putvmorphism(\xpos,\ypos)[\phantom{#3}`#4`#7]{#8}1{#9}
}}

\def\putthreehmorphisms(#1)[#2`#3;#4`#5`#6]#7(#8)#9{{%
\setpos(#1) \settypes(#8)
\if a#9 %
     \vertsize{\tempcounta}{#5}%
     \vertsize{\tempcountb}{#6}%
     \ifnum \tempcounta<\tempcountb \tempcounta=\tempcountb \fi
\else
     \vertsize{\tempcounta}{#4}%
     \vertsize{\tempcountb}{#5}%
     \ifnum \tempcounta<\tempcountb \tempcounta=\tempcountb \fi
\fi
\advance \tempcounta by 60
\puthmorphism(\xpos,\ypos)[#2`#3`#5]{#7}{\arrowtypeb}{#9}
\advance\ypos by \tempcounta
\puthmorphism(\xpos,\ypos)[\phantom{#2}`\phantom{#3}`#4]{#7}{\arrowtypea}{#9}
\advance\ypos by -\tempcounta \advance\ypos by -\tempcounta
\puthmorphism(\xpos,\ypos)[\phantom{#2}`\phantom{#3}`#6]{#7}{\arrowtypec}{#9}
}}

\def\putarc(#1,#2)[#3`#4`#5]#6#7#8{{%
\xpos #1
\ypos #2
\width #6
\arrowlength #6
\putbox(\xpos,\ypos){#3\vphantom{#4}}%
{\advance \xpos by\arrowlength
\putbox(\xpos,\ypos){\vphantom{#3}#4}}%
\horsize{\tempcounta}{#3}%
\horsize{\tempcountb}{#4}%
\divide \tempcounta by2
\divide \tempcountb by2
\advance \tempcounta by30
\advance \tempcountb by30
\advance \xpos by\tempcounta
\advance \arrowlength by-\tempcounta
\advance \arrowlength by-\tempcountb
\halflength=\arrowlength \divide\halflength by 2
\divide\arrowlength by 5
\put(\xpos,\ypos){\bezier{\arrowlength}(0,0)(50,50)(\halflength,50)}
\ifnum #7=-1 \put(\xpos,\ypos){\vector(-3,-2)0} \fi
\advance\xpos by \halflength
\put(\xpos,\ypos){\xpos=\halflength \advance\xpos by -50
   \bezier{\arrowlength}(0,50)(\xpos,50)(\halflength,0)}
\ifnum #7=1 {\advance \xpos by
   \halflength \put(\xpos,\ypos){\vector(3,-2)0}} \fi
\advance\ypos by 50
\vertsize{\tempcounta}{#5}%
\divide\tempcounta by2
\advance \tempcounta by20
\if a#8 %
   \advance \ypos by\tempcounta
   \putbox(\xpos,\ypos){#5}%
\else
   \advance \ypos by-\tempcounta
   \putbox(\xpos,\ypos){#5}%
\fi
}}

\makeatother

\usepackage{amsthm}
\usepackage{dsfont}
\usepackage{stmaryrd}
\newtheorem{theorem}{Theorem}[section]
\newtheorem{lemma}[theorem]{Lemma}
\newtheorem{corollary}[theorem]{Corollary}

\newtheorem{proposition}[theorem]{Proposition}

\begin{document}

\sloppy

%commands
\newcommand{\nl}{\hspace{2cm}\\ }

\def\nec{\Box}
\def\pos{\Diamond}
\def\diam{{\tiny\Diamond}}

\def\lc{\lceil}
\def\rc{\rceil}
\def\lf{\lfloor}
\def\rf{\rfloor}
\def\lk{\langle}
\def\rk{\rangle}
\def\blk{\dot{\langle\!\!\langle}}
\def\brk{\dot{\rangle\!\!\rangle}}

\newcommand{\pa}{\parallel}
\newcommand{\lra}{\longrightarrow}
\newcommand{\hra}{\hookrightarrow}
\newcommand{\hla}{\hookleftarrow}
\newcommand{\ra}{\rightarrow}
\newcommand{\la}{\leftarrow}
\newcommand{\lla}{\longleftarrow}
\newcommand{\da}{\downarrow}
\newcommand{\ua}{\uparrow}
\newcommand{\dA}{\downarrow\!\!\!^\bullet}
\newcommand{\uA}{\uparrow\!\!\!_\bullet}
\newcommand{\Da}{\Downarrow}
\newcommand{\DA}{\Downarrow\!\!\!^\bullet}
\newcommand{\UA}{\Uparrow\!\!\!_\bullet}
\newcommand{\Ua}{\Uparrow}
\newcommand{\Lra}{\Longrightarrow}
\newcommand{\Ra}{\Rightarrow}
\newcommand{\Lla}{\Longleftarrow}
\newcommand{\La}{\Leftarrow}
\newcommand{\nperp}{\perp\!\!\!\!\!\setminus\;\;}
\newcommand{\pq}{\preceq}

\newcommand{\lms}{\longmapsto}
\newcommand{\ms}{\mapsto}
\newcommand{\subseteqnot}{\subseteq\hskip-4 mm_\not\hskip3 mm}

\def\o{{\omega}}

\def\bA{{\bf A}}
\def\bEM{{\bf EM}}
\def\bM{{\bf M}}
\def\bN{{\bf N}}
\def\bF{{\bf F}}
\def\bC{{\bf C}}
\def\bI{{\bf I}}
\def\bK{{\bf K}}
\def\bL{{\bf L}}
\def\bT{{\bf T}}
\def\bS{{\bf S}}
\def\bD{{\bf D}}
\def\bB{{\bf B}}
\def\bW{{\bf W}}
\def\bP{{\bf P}}
\def\bX{{\bf X}}
\def\bY{{\bf Y}}
\def\ba{{\bf a}}
\def\bb{{\bf b}}
\def\bc{{\bf c}}
\def\bd{{\bf d}}
\def\bh{{\bf h}}
\def\bi{{\bf i}}
\def\bj{{\bf j}}
\def\bk{{\bf k}}
\def\bm{{\bf m}}
\def\bn{{\bf n}}
\def\bp{{\bf p}}
\def\bq{{\bf q}}
\def\be{{\bf e}}
\def\br{{\bf r}}
\def\bi{{\bf i}}
\def\bs{{\bf s}}
\def\bt{{\bf t}}
\def\jeden{{\bf 1}}
\def\dwa{{\bf 2}}
\def\trzy{{\bf 3}}

\def\cB{{\cal B}}
\def\cA{{\cal A}}
\def\cC{{\cal C}}
\def\cD{{\cal D}}
\def\cE{{\cal E}}
\def\cEM{{\cal EM}}
\def\cF{{\cal F}}
\def\cG{{\cal G}}
\def\cI{{\cal I}}
\def\cJ{{\cal J}}
\def\cK{{\cal K}}
\def\cL{{\cal L}}
\def\cN{{\cal N}}
\def\cM{{\cal M}}
\def\cO{{\cal O}}
\def\cP{{\cal P}}
\def\cQ{{\cal Q}}
\def\cR{{\cal R}}
\def\cS{{\cal S}}
\def\cT{{\cal T}}
\def\cU{{\cal U}}
\def\cV{{\cal V}}
\def\cW{{\cal W}}
\def\cX{{\cal X}}
\def\cY{{\cal Y}}

%categories

%of functors and monads
\def\Mnd{{\bf Mnd}}
\def\AMnd{{\bf AnMnd}}
\def\An{{\bf An}}
\def\San{{\bf San}}
\def\PMnd{{\bf PolyMnd}}
\def\SanMnd{{\bf SanMnd}}
\def\RiMnd{{\bf RiMnd}}
\def\End{{\bf End}}

%of theories
\def\ET{\bf ET}
\def\RegET{\bf RegET}
\def\RET{\bf RegET}
\def\LrET{\bf LrET}
\def\RiET{\bf RiET}
\def\SregET{\bf SregET}
\def\Cart{\bf Cart}
\def\wCart{\bf wCart}
\def\CartMnd{\bf CartMnd}
\def\wCartMnd{\bf wCartMnd}

%of Lawvere theories
\def\LT{\bf LT}
\def\RegLT{\bf RegLT}
\def\ALT{\bf AnLT}
\def\RiLT{\bf RiLT}

%of Operads
\def\FOp{\bf FOp}
\def\RegOp{\bf RegOp}
\def\SOp{\bf SOp}
\def\RiOp{\bf RiOp}

%various other categories and such
\def\bCat{{{\bf Cat}}}
\def\MonCat{{{\bf MonCat}}}
\def\Mon{{{\bf Mon}}}
\def\Cat{{{\bf Cat}}}

\def\F{\mathds{F}}
\def\S{\mathds{S}}
\def\I{\mathds{I}}
\def\B{\mathds{B}}

%functors
\def\V{\mathds{V}}
\def\W{\mathds{W}}
\def\M{\mathds{M}}
\def\N{\mathds{N}}
\def\R{\mathds{R}}

\def\Op{{\cal O}p}

\def\Vb{\bar{\mathds{V}}}
\def\Wb{\bar{\mathds{W}}}
\def\Sym{{\cal S}ym}

%leftovers
\def\P{{\cal P}}
\def\Q{{\cal Q}}

\pagenumbering{arabic} \setcounter{page}{1}

\title{\bf\Large Monads of regular theories}

\author{ Stanis\l aw Szawiel,\\ Marek Zawadowski\\
%Instytut Matematyki, Uniwersytet Warszawski\\
%ul. S.Banacha 2,\\
%00-913 Warszawa, Poland\\
%s.szawiel@students.mimuw.edu.pl\\
%zawado@mimuw.edu.pl\\
%\date{May 25, 2012\footnote{Corrections May 29, 2012}}
%\date{June 15, 2012}
}

\maketitle
\begin{abstract}  We characterize the category of monads on $Set$ and the category of Lawvere theories that are equivalent to the category of regular equational theories.
\end{abstract}

\tableofcontents

\section{Introduction}

The category of algebras of a (finitary) equational theory can be equivalently described as either a category of models of a Lawvere theory or as a category of algebras of a finitary monad on the category $Set$ or a category of algebras of a (generalized) operad, in which not only permutations but all functions between finite sets act on operations. In fact, the four categories of (finitary) equational theories, Lawvere theories, finitary monads on $Set$ and (generalized) operads are equivalent. These equivalences induce a correspondence between various subcategories.  In \cite{SZ} we have given an intrinsic characterizations of equational theories and Lawvere theories that correspond to the analytic and polynomial monads on $Set$. This was achieved via correspondence with symmetric and rigid operads.

Recall that an equational theory is regular if it has a set of axioms in which each equation contains the same variable on both sides. Thus the theories of monoids and of sup-lattices are regular but the theory of groups is not. The category of regular equational theories and regular morphisms $\RegET$ was defined in \cite{SZ} but the notion of a regular theory was discovered and studied in universal algebra (cf. \cite{Pl}) and will turn fifty soon.

The main objective of this paper is to describe the categories of regular Lawvere theories $\RegLT$, of semi-analytic monads $\SanMnd$, and regular operads $\RegOp$ that correspond to the category of regular theories.

In \cite{SZ} symmetric operads proved to be very useful in describing correspondences of this kind. Full and regular operads (cf. \cite{Tr}, Section \ref{pres_of_cats_of_algs}) play the same role here. The category of full operads can be thought of as the category of monoids in the monoidal category $Set^{\F}$, where $\F$ is the skeleton of the category of finite sets, with the substitution tensor. Thus a monoid in $Set^{\F}$ consists of operations that are equipped with not only a multiplication operator but also an action of the whole category $\F$. As the left Kan extension $Set^{\F}\ra \End$ along the inclusion $\F\ra Set$ is an equivalence of monoidal categories ($\End$ is the category of finitary endofunctors on $Set$) the category of full operads is equivalent to the category of finitary monads on $Set$. Regular operads can be thought of as monoids in the monoidal category $Set^\S$ with substitution tensor, where $\S$ is the skeleton of the category of finite sets and surjections. Thus in regular operads only surjections act on operations.

We identify the essential image of the left Kan extension $Set^{\S}\ra \End$ along inclusion $\S\ra Set$, as the category of semi-analytic functors and semi-analytic natural transformations $\San$. Semi-analytic functors have similar characterizations as analytic ones; cf. Section \ref{pres_of_cats_of_algs}. A finitary \mbox{endofunctor} on $Set$ is semi-analytic iff it preserves pullbacks along monomorphisms. Semi-analytic functors also have presentations via series similar to the series that represent analytic functors.  A natural transformation is semi-analytic iff the naturality squares for monomorphism are pullbacks. The category of semi-analytic functors $\San$ is a monoidal subcategory of $\End$ and the monoids there form the category of semi-analytic monads equivalent to the category of regular theories; cf. Section \ref{sec_mnd}.

The category of regular Lawvere theories is defined very much in the spirit of the definition of the category of analytic Lawvere theories, cf. Section \ref{pres_of_cats_of_algs}. A regular Lawvere theory is a Lawvere theory with nicely behaving isomorphisms and a factorization system consisting of the class of projections and the class of regular morphisms. Regular morphisms in a Lawvere theory are by definitions those morphisms that are right orthogonal to all projections. A regular interpretation of Lawvere theories is an interpretation of Lawvere theories that preserves regular morphisms. We show that the category of regular Lawvere theories is equivalent to the category of regular operads and hence also to the category of semi-analytic monads and to the category of regular equational theories; cf. Section \ref{sec_eqth}.

The category of semi-analytic monads $\SanMnd$ contains subcategories of cartesian $\CartMnd$ and weakly cartesian $\wCartMnd$ monads (i.e. those that preserve and weakly preserve pullbacks, respectively). One may think that such categories of monads can correspond to some natural subcategories of the category of Lawvere theories and of the category of equational theories. In Section \ref{sec_cart} we characterize the subcategories of $Set^\S$ that have as their essential images $\wCartMnd$ and $\CartMnd$. As these characterizations are a bit technical we do not try to rephrase those conditions in terms of either Lawvere theories or equational theories.

In this way we shall describe the second level (level r) of the following diagram
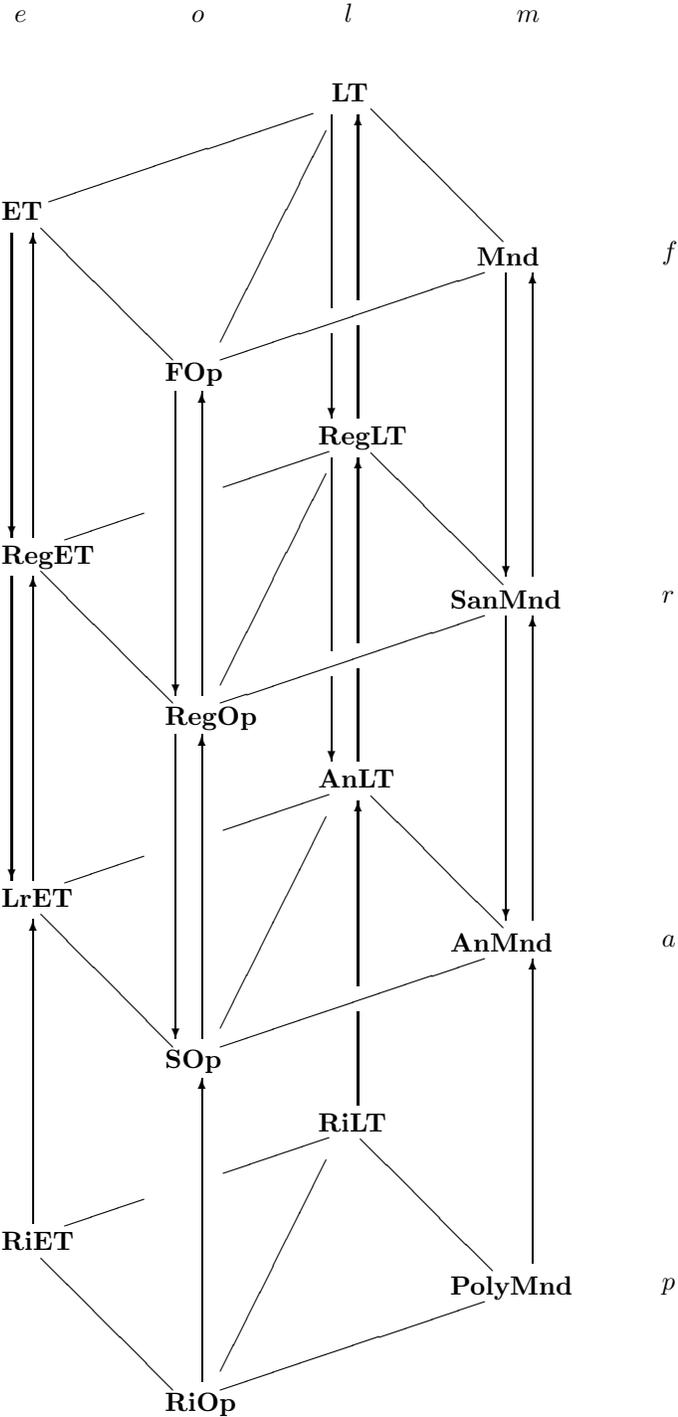
\begin{center} \xext=2800 \yext=5350
\begin{picture}(\xext,\yext)(\xoff,\yoff)
%level p
  %cats
  \put(620,40){$\RiOp$}
  \put(1700,480){$\PMnd$}
  \put(0,650){$\RiET$}
  \put(1200,1100){$\RiLT$}

 %horizontal
 \put(830,120){\line(3,1){1000}}
 \put(840,940){\line(3,1){400}}  \put(540,840){\line(-3,-1){300}}
 \put(650,120){\line(-1,1){500}}
 \put(1860,570){\line(-1,1){500}}
 \put(830,190){\line(1,2){400}}

 %vertical
   \put(120,750){\vector(0,1){1150}}
   \put(760,150){\vector(0,1){1150}}
   \put(2010,600){\vector(0,1){1150}}
   \put(1350,1200){\line(0,1){350}}
   \put(1350,1650){\vector(0,1){700}}

%level a
 %cats
  \put(620,1340){$\SOp$}
  \put(1700,1780){$\AMnd$}
  \put(0,1950){$\LrET$}
  \put(1200,2400){$\ALT$}

 %horizontal
 \put(830,1420){\line(3,1){1000}}
 \put(840,2240){\line(3,1){400}}  \put(540,2140){\line(-3,-1){300}}
 \put(650,1420){\line(-1,1){500}}
 \put(1900,1870){\line(-1,1){500}}
 \put(830,1490){\line(1,2){400}}

   %vertical
   \put(40,3200){\vector(0,-1){1150}} \put(120,2050){\vector(0,1){1150}}
   \put(660,2600){\vector(0,-1){1150}} \put(760,1450){\vector(0,1){1150}}
   \put(1910,3050){\vector(0,-1){1150}} \put(2010,1900){\vector(0,1){1150}}

    \put(1250,3650){\line(0,-1){730}} \put(1250,2820){\vector(0,-1){320}}
    \put(1350,2500){\line(0,1){350}} \put(1350,2950){\vector(0,1){700}}

%level r
 %cats
  \put(620,2640){$\RegOp$}
  \put(1700,3080){$\SanMnd$}
  \put(0,3250){$\RegET$}
  \put(1200,3700){$\RegLT$}

 %horizontal
 \put(830,2720){\line(3,1){1000}}
 \put(840,3540){\line(3,1){400}}  \put(540,3440){\line(-3,-1){300}}
 \put(650,2720){\line(-1,1){500}}
 \put(1900,3170){\line(-1,1){500}}
 \put(830,2790){\line(1,2){400}}

 %vertical
 \put(40,4500){\vector(0,-1){1150}} \put(120,3350){\vector(0,1){1150}}
 \put(660,3900){\vector(0,-1){1150}} \put(760,2750){\vector(0,1){1150}}
 \put(1910,4350){\vector(0,-1){1150}} \put(2010,3200){\vector(0,1){1150}}

 \put(1250,4950){\line(0,-1){730}} \put(1250,4120){\vector(0,-1){320}}
 \put(1350,3800){\line(0,1){350}} \put(1350,4250){\vector(0,1){700}}

%level f
 %cats
  \put(620,3940){$\FOp$}
  \put(1800,4380){$\Mnd$}
  \put(0,4550){$\ET$}
  \put(1250,5000){$\LT$}

 %horizontal
 \put(830,4020){\line(3,1){1000}}
 \put(180,4620){\line(3,1){1000}}
 \put(650,4020){\line(-1,1){500}}
 \put(1900,4470){\line(-1,1){500}}
 \put(830,4090){\line(1,2){400}}
  \put(2500,4400){$f$}
  \put(2500,3100){$r$}
  \put(2500,1800){$a$}
   \put(2500,500){$p$}

\put(50,5300){$e$}
\put(720,5300){$o$}
\put(1300,5300){$l$}
\put(1950,5300){$m$}
\end{picture}
\end{center}
The vertical lines denote adjoint equivalences. Thus up to equivalence there are only four categories in it, one on each level. One equivalent to the category of all finitary monads on $Set$, second  equivalent to the category of all semi-analitic monads on $Set$, third equivalent to the category of all analytic monads on $Set$, and forth equivalent to the category of all polynomial monads on $Set$.  These levels are denoted by letters $f$, $r$, $a$, and $p$, respectively. Thus all four columns of equational theories, Lawvere theories, monads and operads are equivalent. These columns are denoted by letters $e$, $l$, $m$ and $o$, respectively. The vertical functors heading up, in all columns but the column of operads, are inclusions of subcategories. In the column of operads the functors heading up are more like free extensions of the actions. The lower functors are full embeddings and the upper are embeddings that are full on isomorphisms. The vertical functors heading down, the right adjoints to those heading up are monadic. All the squares in the diagram commute up to canonical isomorphisms. The notation concerning categories involved is displayed in the above diagram.

The notation concerning functors is not on the diagram but it is meant to be systematic referring to levels and columns they 'connect'. The horizontal functors are denoted using letters from both columns they connect; the codomain by the script letter, the domain by its subscript, and the level is denoted by superscript. Thus the functor $\AMnd\ra \ALT$ will be denoted by $\cL_m^a$. We usually drop superscripts when it does not lead to confusion. Thus we can write, for example, $\cE_o=\cE_o^p : \RiOp \ra \RiET$. The vertical functors heading up are denoted by the script letter $\cP$ with superscript indicating the column and subscript indicating the level of the codomain. The  vertical functors heading down are denoted by the script letter $\cQ$ with subscript and superscript as with those heading up. Thus we have, for example, functors $\cP=\cP^o=\cP^o_a: \RiOp \ra \SOp$ and $\cQ=\cQ_a=\cQ_a^m: \Mnd \ra \AMnd$. We will also refer to various diagonal morphisms and then we need to extend the notation concerning vertical functors by specifying both the columns of the domain and the codomain. For example, we write $\cP^{ol}_f: \SOp\ra \LT$ to denote one such functor and its right adjoint will be denoted by $\cQ^{lo}_a: \LT\ra \SOp$. In principle this notation will leave the codomain not always uniquely specified but in practice it it sufficient, and in fact, usually much less is needed.

\subsection*{Notation}
$n=\{0,\ldots, n-1\}$, $[n]=\{0,\ldots, n\}$, $(n]=\{1,\ldots, n\}$, $\o$ - denotes the set of natural numbers. The set $X^n$ is interpreted as $X^{(n]}$
when convenient.
The skeletal category equivalent to the category of finite sets $Set_{fin}$ will be denoted by $\mathds{F}$.
We will be assuming that the objects of $\mathds{F}$ are sets $(n]$ for $n\in \o$.  The subcategories of $\mathds{F}$ with
the same objects as $\mathds{F}$ but having as objects bijection, surjections and injections will be denoted by
$\mathds{B}$, $\mathds{S}$, $\mathds{I}$, respectively. When $S_n$ acts on a set $A_n$ on the right and on the set $B_n$ on the left, the set $A\otimes_nB$
is the usual tensor product of $S_n$-sets.

%\newpage
\section{Presentations of categories of algebras}\label{pres_of_cats_of_algs}
\subsection*{Equational theories}
The category of regular theories, denoted by $\RegET$,  was introduced in \cite{SZ}. It is a subcategory of all equational theories, denoted $\ET$.

\subsection*{Operads}
The symmetric operads provide yet another way of presenting models of an equational theory. This kind of presentation is usually very convenient, however the models defined by such operads are more specific than models of arbitrary equational theories. For example, if $\cO$ is a symmetric operad then the free algebra functor $Set\ra Alg(\cO)$ preserves weak wide pullbacks. Below we extend the definition of an operad so that it captures all the equational theories but still keeps the operadic flavor. The main difference is that instead of having just symmetric groups acting on sets of operations we have actions of the morphisms of the whole skeleton of the category of finite sets $\F$. Symmetric operads can be thought of as monoids for the substitution tensor on the category $Set^\B$. Similarly $\F$-operads can be thought of as monoids for the substitution tensor on the category $Set^\F$. By $\End$ we denote the category of finitary endofunctors of $Set$. It is a strict monoidal category with the tensor being composition. The substitution tensor on $Set^\F$ makes the equivalence of categories $Set^\F \ra \End$, given by left Kan extension, a strong monoidal equivalence. This immediately shows that the category of $\F$-operads is equivalent to the category of finitary monads on $Set$. The definition below was first explicitly spelled out in \cite{Tr}.

A {\em full operad} (or $\F$-operad) $\cO$ consists a family of sets $\cO_n$, for $n\in \o$, a unit element
$\iota\in \cO_1$ for any $k,n,n_1,\ldots, n_k \in \o$ with $n={\sum_{i=1}^k n_i}$ a multiplication operation
\[  \ast : \cO_{n_1}\times\ldots \times \cO_{n_k} \times  \cO_k \lra \cO_n\]
\[ (b_1,\ldots, b_k,a) \mapsto \lk b_1,\ldots, b_k\rk\ast a \]
a left action of the morphisms on operations
\[  \cdot : \F(n,m)\times \cO_n \lra \cO_m \]
for $n,m\in \o$, such that the multiplication is associative with unit $\iota$ and compatible with the category action, i.e.
for $a\in \cO_n$
\[ \lk a \rk \ast \iota = a = \lk \iota, \ldots, \iota \rk \ast a; \;\;\;\; 1_n\cdot a = a\]
for $a\in\cO_n$, $b_i\in O_{k_i}$,  $c_{ij}\in \cO_{m_{ij}}$, for $i\in (n]$, $j\in (k_i]$ we have
\[ \lk c_{1,1},\ldots, c_{1,k_1},\ldots,c_{n,1},\ldots, c_{n,k_n} \rk  \ast (\lk b_1,\ldots, b_n\rk \ast a) =\]
\[=\lk  \lk c_{1,1},\ldots, c_{1,k_1}\rk \ast b_1,\ldots, \lk c_{n,1},\ldots, c_{n,k_n}\rk \ast b_n \rk \ast a \]
and for $\phi_\in \F(n,m)$,  $a\in \cO_n$ $\psi_i\in \F(k_i,l_i)$,  $b_i\in \cO_{k_i}$ for $i\in (m]$, we have
\[ \lk \psi_1\cdot b_1, \ldots \psi_n\cdot b_n \rk \ast (\phi\cdot a) =
(\lk \psi_1\ldots , \psi_n\rk \star \phi)\cdot (\lk b_{\phi(1)} , \ldots, b_{\phi(n)}\rk \ast a) \]
in $\cO_l$, where
\[ \lk \psi_1\ldots , \psi_n\rk\star \phi : (k]=\sum_{j=1}^n (k_{\phi(j)}]\lra \sum_{i=1}^m (l_i]=(l] \]
is a function such that $j$-th summand of the domain $(k_{\phi(j)}]$ is sent to the $\phi(j)$-th summand  $(l_{\phi(j)}]$ of the codomain by the function $\psi_{\phi(j)}$, i.e.
\[ \lk j,r\rk\mapsto \lk\phi(j),\psi_{\phi(j)}(r) \rk \]
for $j\in (n]$ and $r\in (k_{\phi(j)}]$. This definition refers to the obvious lexicographic order on  both $(k]=\sum_{j=1}^n (k_{\phi(j)}]$ and  $(l]=\sum_{i=1}^m (l_i]$. It is an extension of the multiplication operation in the operad of symmetries.

\begin{center} \xext=2000 \yext=450
\begin{picture}(\xext,\yext)(\xoff,\yoff)
\setsqparms[1`-1`-1`1;2000`400]
 \putsquare(0,0)[\sum_{j=1}^n (k_{\phi(j)}]`\sum_{i=1}^m (l_i]`(k_{\phi(j)}]`(l_{\phi(j)};\lk \psi_1\ldots , \psi_n\rk\star \phi`\kappa_j`\lambda_{\phi(j)}`\psi_{\phi(j)}]
 \end{picture}
\end{center}

{\bf Remark} The operation $\star$ is defined on a family of function indexed by $k,l\in\o$ between sets:
\begin{center}
\xext=2200 \yext=150
\begin{picture}(\xext,\yext)(\xoff,\yoff)
\putmorphism(600,50)(1,0)[\coprod_{n,m\in\o,\; \phi\in \F(n,m)}\;\;\coprod_{k_i, l_i\in\o,\;
k=\sum_{j=1}^n k_{\phi(j)},\; l=\sum_{i=1}^m l_i } \;\;\prod_{i=1}^m \F(k_i,l_i)`\F(k,l)`\star]{2200}{1}a
\end{picture}
\end{center}

A {\em morphism of full operads} $f: \cO \ra \cO'$ is a function that respects arities of operations, unit, compositions, and the actions of functions from $\F$.

The operation $\star$, when applied to morphisms in $\S$, returns a morphism in $\S$. Therefore these definitions make sense if we restrict morphisms to surjections i.e. morphisms in $\S$\footnote{Clearly, these definitions make also sense if we restrict morphisms in $\F$ to injections i.e. morphisms in $\I$. But we will study such operads elsewhere.}. In this way, we obtain the notion of a {\em regular operad}, a {\em morphism of regular operads} and the whole category of regular operads denoted $\RegOp$.

We have functors
\begin{center}
\xext=1800 \yext=100
\begin{picture}(\xext,\yext)(\xoff,\yoff)
\putmorphism(0,20)(1,0)[\FOp`\RegOp`\cQ_r^o]{800}{1}a
\putmorphism(800,20)(1,0)[\phantom{\RegOp}`\SOp`\cQ_a^o]{800}{1}a
\end{picture}
\end{center}
'restricting actions' along the inclusions $\B\ra \S\ra \F$. They have left adjoints:
\begin{center}
\xext=1800 \yext=100
\begin{picture}(\xext,\yext)(\xoff,\yoff)
\putmorphism(0,20)(1,0)[\FOp`\RegOp`\cP_f^o]{800}{-1}a
\putmorphism(800,20)(1,0)[\phantom{\RegOp}`\SOp`\cP_r^o]{800}{-1}a
\end{picture}
\end{center}
We sketch the definitions of those functors below.

For a symmetric operad $\cO$ the regular operad $\cP_r^o(\cO)$ has
\[ \coprod_{k\in \o} \S(k,n)\otimes_k \cO_k \]
as the set of $n$-ary operations. The unit of $\cP_r^o(\cO)$ is $$[1,\iota]\in  \S(1,1)\otimes_1 \cO_1 \subset \coprod_{k\in \o} \S(k,1)\otimes_k \cO_k$$
The action of the category $\S$ in $\cP_r^o(\cO)$
\[ \cdot :\S(n,m)\times \coprod_{k\in \o} \S(k,n)\otimes_k \cO_k \lra \coprod_{k\in \o} \S(k,n)\otimes_k \cO_k \]
is
\[  \phi \cdot [f,a]  = [\phi\circ f,a] \]
where $\phi:(n]\ra (m]$, $f:(k]\ra (n]$ are surjection and $a\in \cO_k$.
The composition in $\cP_r^o(\cO)$
\[ \coprod_{l_1\in \o} \S(l_1,n_1)\otimes_{l_1} \cO_{l_1}\times\ldots \times \coprod_{l_m\in \o} \S(l_m,n_m)\otimes_{l_m} \cO_{l_m}  \times \coprod_{k\in \o} \S(k,m)\otimes_k \cO_k \lra \coprod_{k\in \o} \S(k,n)\otimes_k \cO_k \]
is given by
\[ \lk [\psi_1,b_1],\ldots , [\psi_m,b_m],[\phi,b] \rk \mapsto  [ \lk \psi_1,\ldots, \psi_m\rk\star\phi, \lk b_{\phi(1)},\ldots,b_{\phi(k)} \rk \ast a ]  \]
where $\varphi\in\S(k,m)$, $a\in \cO_{k}$, $\sum_{i=1}^m n_i =n$, $\psi_i\in\S(l_i,n_i)$, $b_i\in \cO_{l_i}$, for  $i\in (m]$. The definition of the functor $\cP_a^o$ on morphisms is left for the reader.

For a regular operad $\cO$ the full operad $\cP_f^o(\cO)$ has
\[ \coprod_{k\in \o} \I(k,n)\otimes_k \cO_k \]
as the set of $n$-ary operations.
The action of the category $\F$ in $\cP_f^o(\cO)$
\[ \cdot :\F(n,m)\times \coprod_{k\in \o} \I(k,n)\otimes_k \cO_k \lra \coprod_{k\in \o} \I(k,n)\otimes_k \cO_k \]
is
\[  \phi \cdot [f,a]  = [f', \phi'\cdot a] \]
where $\phi:(n]\ra (m]$ is a function, $f:(k]\ra (n]$ is an injection, $a\in \cO_k$ and  $\phi'$, $f'$ is the epi-mono factorization of $\phi\circ f$:
\begin{center} \xext=800 \yext=500
\begin{picture}(\xext,\yext)(\xoff,\yoff)
\setsqparms[2`3`1`2;800`400]
 \putsquare(0,50)[(k]`(n]`{(n']}`{(m]};f`\phi'`\phi`f']
 \end{picture}
\end{center}
The rest of the definition of  $\cP_f^o(\cO)$ is similar to the remaining part of the definition of  $\cP_r^o(\cO)$ given above.

\subsection*{Lawvere theories}
The category of Lawvere theories will be denoted by $\LT$, see \cite{Lw}, \cite{KR}, \cite{SZ} for details.  $\F^{op}$ is the initial Lawvere theory with the obvious projections. Let $\bT$ be any Lawvere theory.  By $\pi : \F^{op}\ra \bT$ we  denote the unique morphism from $\F^{op}$ to $\bT$.  The class of {\em projections} in $\bT$ is the closure under isomorphisms of the image of the injections in $\F$ under $\pi$.  A morphism $r$ in $T$ is {\em regular} iff $r$ is right orthogonal to all projection morphisms in $\bT$. By a factorization system we mean a factorization system in the sense of [FK], see [CJKP] sec 2.8.

$Aut(n)$ is the set of automorphisms of $n$ in $\bT$. Recall from \cite{SZ}, that a Lawvere theory has {\em simple automorphisms} if the canonical function
\[ \rho_n : S_n \times Aut(1)^n \lra Aut(n) \]
such that
\[ (\sigma, a_1,\ldots, a_n) \mapsto a_1\times  \ldots \times a_n \circ \pi_\sigma \]
is a bijection, for $n\in\o$.

A Lawvere theory $\bT$ is a {\em regular Lawvere theory} iff the projection morphisms and the regular morphisms form a factorization system and $\bT$ has simple automorphisms. A {\em regular interpretation} of Lawvere theories is a morphism of Lawvere theories that preserves regular morphisms. Thus we have a non-full subcategory of $\LT$ of regular Lawvere theories with regular interpretations $\RegLT$. The theory $\F^{op}$ is regular.

We have inclusion functors
\begin{center}
\xext=1800 \yext=100
\begin{picture}(\xext,\yext)(\xoff,\yoff)
\putmorphism(0,20)(1,0)[\ALT`\RegLT`\cP_r^l]{800}{1}a
\putmorphism(800,20)(1,0)[\phantom{\RegLT}`\LT`\cP_f^l]{800}{1}a
\end{picture}
\end{center}

\subsection*{Monads}\label{monads} We introduce here a notion of a semi-analytic monad that is broader then that of an analytic monad but it still retains some combinatorial flavor.

Recall that an analytic (endo)functor on $Set$ can be defined by any of the following conditions
\begin{enumerate}
  \item finitary functor preserving weak wild pullbacks;
  \item Kan extension of a functor from $\B$ to $Set$;
  \item functor (isomorphic to one) having an analytic presentation $\sum_{n\in \o} X^n\otimes_n A_n$, where the $n$-coefficient $A_n$ is a (left) $S_n$-set for $n\in\o$.
\end{enumerate}

Similarly we shall define a {\em semi-analytic} (endo)functor on $Set$ as a functor satisfying any of three equivalent conditions (see Proposition \ref{hat vs kan}, Theorem \ref{im hat})
\begin{enumerate}
  \item finitary functor preserving pullbacks along monomorphisms;
  \item Kan extension of a functor from $\S$ to $Set$ along the inclusion functor $\bi_\S:\S\ra Set$;
  \item functor (isomorphic to one) having a semi-analytic presentation $\sum_{n\in \o} \left[\begin{array}{c} X \\ n  \end{array}\right] \otimes_n A_n$, where the category $\S$ acts on coefficients $A_n$ on the left (see below).
\end{enumerate}
The first two characterizations of the semi-analytic functors are clear. We shall describe the third one and show that it is equivalent to the other two. For the time being, the first definition of a semi-analytic functor is the official one. A natural transformation $\phi: F\ra G$ is {\em semi-cartesian} iff the naturality squares for monomorphisms are pullbacks. The category of semi-analytic functors with the semi-cartesian natural transformations will be denoted by  $\San$.
\vskip 2mm
{\em Examples}
\begin{enumerate}
  \item The functor $\cP_{\leq n}:Set \lra Set$  associating to a set $X$ the set of subsets of $X$ with at most $n$-elements is not analytic, if $n>2$, as it can be easily seen that it does not preserve weak pullbacks. However, it preserves pullbacks along monos and hence it is semi-analytic.
  \item If $U$ is a set, $n\in \o$ then the functor $(-)_{\leq n}^U:Set \ra Set$, associating to a set $X$ the set of functions from $U$ to $X$ with an at most $n$-element image, is not analytic, if $|U|>n>2$. Again it can be easily seen that it does not preserve weak pullbacks. However, it is semi-analytic.
  \item We will see later that the functor part of any monad on $Set$ that comes from a regular equational theory is semi-analytic.
\end{enumerate}

For a set $X$ and $n\in \o$ the set
\[ \left[\begin{array}{c} X \\ n  \end{array}\right] \]
denotes the set of monomorphisms from $(n]$ to $X$. If $X$ has less than $n$ elements, this set is empty and if $X$ is a finite set, then it has $\left(\begin{array}{c} |X| \\ n  \end{array}\right)\cdot k!$ elements. The notation is in analogy with the notation
\[ \left(\begin{array}{c} X \\ n  \end{array}\right) \]
denoting the set of $n$-element subsets of $X$.

Clearly, $\left[\begin{array}{c} X \\ n  \end{array}\right]$ is not functorial in $X$ on its own but if we build a series with such sets and coefficients that are related by surjections we do get a functor. To see this let
\[ A: \S \ra Set\]
be a functor, or equivalently a sequence of coefficient sets $\{ A_n \}_{n\in\o}$ on which the category $\S$ acts on the left. As $S_n$ acts on $\left[\begin{array}{c} X \\ n  \end{array}\right]$ on the right, by composition it makes sense to form a set
\[ \left[\begin{array}{c} X \\ n  \end{array}\right] \otimes_n A_n\]
and a whole coproduct
\[ \hat{A}(X)=\sum_{n\in\o}\left[\begin{array}{c} X \\ n  \end{array}\right] \otimes_n A_n \]
The above formula is functorial in $X$, i.e. $\hat{A}$ can be defined on morphisms as follows. Let $f:X\ra Y$ be a function and  $[\vec{x},a]$ an element of $\left[\begin{array}{c} X \\ n  \end{array}\right] \otimes_n A_n$. We take the epi-mono factorization $\alpha$, $\vec{y}$ of $f\circ \vec{x}$
\begin{center} \xext=500 \yext=450
\begin{picture}(\xext,\yext)(\xoff,\yoff)
\setsqparms[2`3`1`2;500`400]
 \putsquare(0,0)[(n]`X`(m]`Y;\vec{x}`\alpha`f`\vec{y}]
 \end{picture}
\end{center}
and we put
\[ \hat{A}(f)([\vec{x},a]) = [\vec{y},A(\alpha)(a)] \]
which is an element in $\left[\begin{array}{c} Y \\ m  \end{array}\right] \otimes_m A_m$, i.e. in $\hat{A}(Y)$. As the factorization is unique up to an isomorphism, $\hat{A}(f)$ is well defined.

If $\tau: A\ra B$ is a natural transformation of functors it induces a natural transformation of functors
\[ \hat{\tau} : \hat{A}\lra \hat{B} \]
as follows. For  $[\vec{x},a]$ in $\left[\begin{array}{c} X \\ n  \end{array}\right] \otimes_n A_n$ we put
\[ \hat{\tau}([\vec{x},a])= [\vec{x},\tau(a)] \]

\begin{proposition}\label{hat vs kan}
The functors $\hat{(-)}: Set^\S \lra \End$ is well defined and it is isomorphic to the left Kan extension along the inclusion functor $\bi_\S:\S\ra Set$.
\end{proposition}

{\em Proof.} It is well known (see \cite{CWM}) that the Kan extensions can be calculated with coends.  Thus for $A:\S\ra Set$ and a function $f:X\ra Y$
we have
\begin{center} \xext=2500 \yext=200
\begin{picture}(\xext,\yext)(\xoff,\yoff)

   \putmorphism(0,50)(1,0)[Lan_{\bi_\S}(A)(X)=\int^{n\in\cS}X^n\times A_n`\int^{n\in\cS}Y^n\times A_n=Lan_{\bi_\S}(A)(Y)`Lan_{\bi_\S}(A)(f)]{2500}{1}a
\end{picture}
\end{center}
\[ \|\vec{x},a\| \mapsto \|f\circ\vec{x},a\| \]
where $\|\vec{x},a\|$ is the equivalence class of the equivalence relation $\approx$ on $\sum_{n\in\o} X^n\times A_n$ generated by the relation $\sim$ such that
\[ \lk \vec{y}\circ \phi, a\rk  \sim  \lk \vec{y}, \phi\cdot a\rk \]
for $\vec{y}:(m]\ra X$, $\phi : (n]\ra (m]\in \S$ and $a\in A_n$.  Let us call an representant $\lk\vec{x},a\rk$ of a class $\|\vec{x},a\|$ {\em minimal} iff $\vec{x}$ is injective. Any element $\lk\vec{x},a\rk$ is $\sim$-related to one of form $\lk \vec{y}, f\cdot a\rk$ where $f$, $\vec{y}$ is the surjective-injective factorization of $\vec{x}$, i.e. any element is $\sim$-related to a minimal one. It is easy to see that any two minimal representatives of $\approx$-equivalence class are $\sim$-related.

We define an isomorphism
\[ \kappa_A : \hat{A} \lra Lan_{\bi_\S}(A) \]
so that on a set $X$
\[ (\kappa_A)_X : \hat{A}(X) \lra Lan_{\bi_\S}(A)(X) \]
we put
\[ [ \vec{x},a] \mapsto \| \vec{x},a\| \]
where $\vec{x} :(n]\ra X$ and $a\in A_n$. Taking the minimal elements in the equivalence class is the inverse function. Thus $(\kappa_A)_X$ is a bijection.
It is easy to see that $\kappa_A$ defined this way is natural in $X$, i.e. it is a natural isomorphism.

For a natural transformation $\tau : A\ra B$ in $Set^\S$ the square
\begin{center} \xext=600 \yext=650
\begin{picture}(\xext,\yext)(\xoff,\yoff)
\setsqparms[2`3`1`2;600`500]
 \putsquare(0,50)[\hat{A}`Lan_{\bi_\S}(A)`\hat{B}`Lan_{\bi_\S}(B);\tau_A`\hat{\tau}`Lan_{\bi_\S}(\tau)`\tau_B]
 \end{picture}
\end{center}
obviously commutes, as both compositions of an element $[\vec{x},a]\in \hat{A}(X)$ are equal to  $\|\vec{x},\tau(a)\|$. Thus
 $$\kappa : \hat{(-)} \lra Lan_{\bi_\S} $$
is also a natural isomorphism, as required. $\boxempty$

\begin{theorem} \label{im hat}
The functor $\hat{(-)}: Set^\S \lra \End$ is faithful, full on isomorphisms, and its essential image is the category of semi-analytic functors $\San$.
\end{theorem}

The theorem will be proved via a series of lemmas.

\begin{lemma} \label{im hat1}
The essential image of the functor $\hat{(-)}: Set^\S \lra \End$ is contained in the category $\San$.
\end{lemma}

{\em Proof.} First we check that for $A\in Set^\S$, $\hat{A}$ preserves pullbacks along monos. Let
\begin{center} \xext=500 \yext=500
\begin{picture}(\xext,\yext)(\xoff,\yoff)
\setsqparms[1`1`1`1;500`400]
 \putsquare(0,0)[P`Y`Z`X;f'`\beta'`\beta`f]
 \end{picture}
\end{center}
be a pullback in $Set$. We need to show that the square
\begin{center} \xext=2000 \yext=1200
\begin{picture}(\xext,\yext)(\xoff,\yoff)
\setsqparms[1`1`1`1;2000`800]
 \putsquare(0,200)[\sum_{n\in\o}\left[\begin{array}{c} P \\ n  \end{array}\right] \otimes_n A_n`
 \sum_{n\in\o}\left[\begin{array}{c} Y \\ n  \end{array}\right] \otimes_n A_n`
 \sum_{n\in\o}\left[\begin{array}{c} Z \\ n  \end{array}\right] \otimes_n A_n`
 {\sum_{n\in\o}\left[\begin{array}{c} X \\ n  \end{array}\right] \otimes_n A_n};\hat{A}(f')`\hat{A}(\beta')`\hat{A}(\beta)`\hat{A}(f)]
 \end{picture}
\end{center}
is also a pullback in $Set$. Let ${[\vec{y},a]}\in \left[\begin{array}{c} Y \\ n  \end{array}\right] \otimes_n A_n$ and
${[\vec{z},b]}\in \left[\begin{array}{c} Z \\ m  \end{array}\right] \otimes_m A_m$ be such that
\[ \hat{A}(f)([\vec{z},b]) = \hat{A}(\beta)([\vec{y},a]) \]
Let $\alpha:(m]\ra (n']$ and $\vec{z}':(n]\ra X$ be a surjection-injection factorization of the function $f\circ \vec{z}:(m]\ra X$.
Then there is a $\sigma\in S_n$ such that
\[ \vec{z}\circ \sigma = \beta\circ \vec{y},\hskip 5mm A(\alpha)(b)=A(\sigma)(a) \]
We have
\[ f\circ \vec{z} = \vec{z}'\circ \alpha= \beta\circ \vec{y}\circ \sigma^{-1}\circ \alpha \]
Hence there is a function $\vec{p}: (m]\ra P$ as in the following diagram
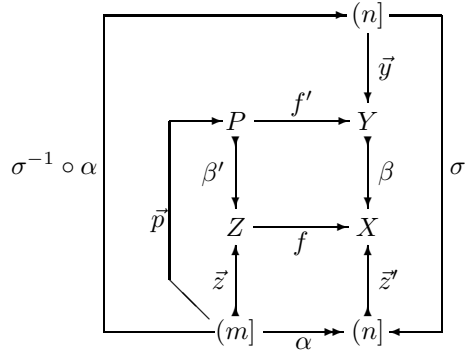
\begin{center} \xext=1200 \yext=1300
\begin{picture}(\xext,\yext)(\xoff,\yoff)
\setsqparms[1`-2`-2`3;500`400]
 \putsquare(500,0)[Z`X`(m]`{(n]};`\vec{z}`\vec{z}'`\alpha]
 \setsqparms[1`2`2`0;500`400]
  \putsquare(500,400)[P`Y`\phantom{Z}`\phantom{X};f'`\beta'`\beta`f]

  \putmorphism(1000,1200)(0,-1)[(n]`\phantom{Y}`\vec{y}]{400}{1}r
  %p
    \put(250,800){\vector(1,0){200}}
    \put(250,800){\line(0,-1){600}}
    \put(250,200){\line(1,-1){150}}
  \put(180,400){$\vec{p}$}
  %\sigma
   \put(1280,0){\vector(-1,0){200}}
    \put(1280,0){\line(0,1){1200}}
    \put(1280,1200){\line(-1,0){200}}
  \put(1310,600){$\sigma$}

    %\sigma^-1..
   \put(0,1200){\vector(1,0){920}}
    \put(0,0){\line(0,1){1200}}
    \put(0,0){\line(1,0){420}}
  \put(-350,600){$\sigma^{-1}\circ \alpha$}
 \end{picture}
\end{center}
such that
\[ \beta'\circ\vec{p}= \vec{z}, \hskip 5mm f'\circ\vec{p}= \vec{y}\circ\sigma^{-1}\circ\alpha \]
The function $\vec{p}$ is an injection since $\vec{z}$ is. Thus $[\vec{p},b]$ is an element of $\left[\begin{array}{c} P \\ m  \end{array}\right] \otimes_m A_m$.
We have
\[ \hat{A}(\beta')([\vec{p},b]) = [\beta'\circ \vec{p},b] = [\vec{z},b] \]
Moreover, on representatives we have
\[ \lk f'\circ\vec{p},b\rk = \lk \vec{y}\circ\sigma^{-1}\circ\alpha,b\rk \sim \]
\[\sim  \lk \vec{y},A(\sigma^{-1}\circ\alpha)(b)\rk = \lk \vec{y},A(\sigma^{-1})( A(\alpha)(b))\rk =\]
\[ =\lk \vec{y},A(\sigma^{-1})( A(\sigma)(a))\rk = \lk \vec{y},a\rk\]
and hence
\[ \hat{A}(f')([\vec{p},b]) = [f'\circ\vec{p},b]= [\vec{y},a] \]
Thus $\hat{A}$ preserves pullbacks along monos  indeed.

Now let $\tau:A\ra B$ be a morphism in $Set^\S$. We shall show that $\hat{\tau}:\hat{A}\ra \hat{B}$ is semi-cartesian.
Let $f:X\ra Y$ be an injection in $Set$. We need to show that the square
\begin{center} \xext=500 \yext=550
\begin{picture}(\xext,\yext)(\xoff,\yoff)
\setsqparms[1`1`1`1;500`400]
 \putsquare(0,50)[\hat{A}(X)`\hat{B}(X)`\hat{A}(Y)`\hat{B}(X);\hat{\tau}_X`\hat{A}(f)`\hat{B}(f)`\hat{\tau}_Y]
 \end{picture}
\end{center}
is a pullback. Let $[\vec{y},a]\in \hat{A}(Y)$ and $[\vec{x},b]\in \hat{B}(X)$ be such that
\[  [\vec{y},\tau_n(a)]=\hat{\tau}_Y([\vec{y},a]) = \hat{B}(f)([\vec{x},b]) =[f\circ\vec{x},b] \]
where $a\in A_n$. Note that as $f$ is an injection, so is $f\circ\vec{x}$ and hence $\lk f\circ\vec{x},b\rk$ represents an element in $\hat{B}(Y)$.
The above equality means that $b\in B_n$ and there is a $\sigma\in S_n$ such that
\[ \vec{y}\circ \sigma = f\circ \vec{x},\hskip 5mm \tau_n(a) = B(\sigma)(b) \]
We have an element $[\vec{x}\circ \sigma^{-1},a]$ such that
\[ \hat{A}(f)([\vec{x}\circ \sigma^{-1},a])= [f\circ \vec{x}\circ \sigma^{-1},a] = [\vec{y},a] \]
On representatives we have
\[ \lk \vec{x},\tau_n(a) \rk  \sim \lk \vec{x},\hat{B}(\sigma^{-1}(\tau_n(a)) \rk = \lk \vec{x},b \rk \]
and hence we have
\[ \hat{\tau}_X([\vec{x}\circ \sigma^{-1},a])=  [\vec{x},\tau_n(a) ] = [ \vec{x},b] \]
as well. This means that $\hat{\tau}$ is semi-cartesian indeed. $\boxempty$

\begin{lemma} \label{im hat2}
The functor $\hat{(-)}: Set^\S \lra \End$ is faithful and full on semi-cartesian morphisms. In particular it is full on isomorphisms.
\end{lemma}

{\em Proof.} One can easily verify that if two natural transformations $\hat{\tau},\hat{\tau}' :\hat{A}\ra \hat{B}$ agree on elements of the form $[1_{(m]},a]$ for $m\in \o$ and $a\in \hat{A}(m]$ then natural transformation $\tau,\tau' :A\ra B$ are equal. Thus $\hat{(-)}$ is faithful.

To show that $\hat{(-)}$ is full, let us fix two functors $A,B\in Set^\S$ and a semi-cartesian natural transformation $\psi :\hat{A}\ra \hat{B}$. Let $m\in \o$ and $a\in A_m$. We shall define $\tau_m(a)$. We have $\psi_{(m]}([1_{(m]},a]) = [f,b]\in \hat{B}(m]$ for some injection $f:(k]\ra (m]$ and $b\in B_k$.
We claim that $k=m$ and $f$ is a bijection. Suppose to the contrary, that $k<m$. We have that $\hat{B}(f)([1_{(k]},b])=[f,b]$. Since $f$ is an injection and $\psi$ is semi-cartesian, the square
\begin{center} \xext=800 \yext=450
\begin{picture}(\xext,\yext)(\xoff,\yoff)
\setsqparms[1`-1`-1`1;800`400]
 \putsquare(0,0)[\hat{A}(m]`\hat{B}(m]`\hat{A}(k]`{\hat{B}(k]};\psi_{(m]}`\hat{A}(f)`\hat{B}(f)`\psi_{(k]}]
 \end{picture}
\end{center}
is a pullback and hence there is an element $[g,c]\in \hat{A}(k]$ such that
\[ \hat{A}(f)([g,c])=[1_{(m]},a], \hskip 5mm \psi_{(m]} ([g,c])=[1_{(k]},b] \]
The first equality implies that the proper mono $f\circ g$ is epi (as its codomain is $(m]$ and its mono part is $1_{(m]}$). This is a contradiction. Hence $f$ is a bijection and we can apply the functor $B$ to it. We put
\[ \tau_m(a)= B(f)(b) \]
Thus we have defined a transformation $\tau:A\ra B$.

To show that $\tau$ is natural, let us fix a surjection $\beta:(m]\ra (n]$ and $a\in A_m$. Then using definitions of $\hat{A}$, $\hat{B}$ and the naturality of $\psi$ we have
\[ [1_{(n]},B(f)(\tau_m(a))] = \hat{B}(f)([1_{(m]},\tau_m(a)]) = \]
\[ =\hat{B}(f)\circ \psi_{(m]}([1_{(m]},a]) = \psi_{(n]}\circ \hat{A}(f)([1_{(m]},a])=\]
\[ = \psi_{(n]}([1_{(n]},\hat{A}(f)(a))]) = [1_{(n]},\tau_n(\hat{A}(f)(a))] \]
But this means that
\[ B(f)(\tau_m(a))= \tau_n(\hat{A}(f)(a)) \]
Since $a$ and $f$ were arbitrary, $\tau$ is natural.

Finally, we show that $\hat{\tau}=\psi$. Let us fix a set $X$ and $[\vec{x}:(k]\ra X,a\in A_k]\in \hat{A}(X)$. Note that $\hat{A}(\vec{x})([1_{(k]},a])=[\vec{x},a]$. Using the naturality of $\psi$ and $\hat{\tau}$ on $\vec{x}$ and the definition of $\tau$ we have
\[ \psi_X([\vec{x},a]) =\psi_X\circ\hat{A}(\vec{x})([1_{(k]},a]) = \]
\[ =\hat{B}(\vec{x})\circ \psi_{(k]}([1_{(k]},a]) = \hat{B}(\vec{x})([1_{(k]},\tau_k(a)])=\]
\[ =\hat{B}(\vec{x})\circ \hat{\tau}_{(k]}(([1_{(k]},a]) =\hat{\tau}_{X}\circ\hat{A}(\vec{x})(([1_{(k]},a])=\]
\[ = \hat{\tau}_{X}([\vec{x},a]) \]
Since $X$ and $[\vec{x},a]$ were arbitrary $\hat{\tau}=\psi$.
$\boxempty$

\begin{lemma} \label{im hat3}
Each semi-cartesian functor is in the essential image of $\hat{(-)}: Set^\S \lra \End$.
\end{lemma}

{\em Proof.} Let us fix a semi-cartesian functor. We define a functor $A:\S\ra Set$ and a natural isomorphism $\tau: \hat{A}\ra F$.

Put
\[ A(n] =F(n] -\bigcup_{(m] \stackrel{f}{\ra}(n]\in \I,\, m<n} F(f)(A(m])  \]
Note that the sum over the empty index set is empty. For $\alpha:[n) \ra [m)$ in $\S$ we put
\[ A(\alpha)=F(\alpha)_{\lceil A[n)} : A[n)\lra A[m) \]
We need to show that $A(\alpha)$ is well defined, i.e. for $a\in A(n]$ we have $F(\alpha)\in A(m]$. Clearly, $F(\alpha)\in F(m]$.  Suppose to the contrary, that $F(\alpha)(a)\not\in A(m]$. Then there is a proper mono $f:(k]\ra (m]$ and $b\in A(k]$ such that $F(\alpha)(a)= F(f)(b)$. We form a pullback of a surjection $\alpha$ and a proper mono $f$
\begin{center} \xext=600 \yext=500
\begin{picture}(\xext,\yext)(\xoff,\yoff)
\setsqparms[3`2`2`3;600`400]
 \putsquare(0,50)[(p]`(k]`(n]`{(m]};\alpha'`f'`f`\alpha]
 \end{picture}
\end{center}
Thus $f'$ is again a proper mono. $F$ preserves this pullback, i.e. the square
\begin{center} \xext=600 \yext=550
\begin{picture}(\xext,\yext)(\xoff,\yoff)
\setsqparms[3`2`2`3;600`400]
 \putsquare(0,50)[F(p]`F(k]`F(n]`{F(m]};F(\alpha')`F(f')`F(f)`F(\alpha)]
 \end{picture}
\end{center}
is again a pullback. Thus there is an element $c\in F(p]$ such that
\[ f(\alpha')(c)= b,\hskip 5mm  f(f')(c)=a  \]
The later equality means that $a\not\in A(n]$  contrary to the supposition. Thus $A:\S\ra Set$ is a well defined functor.

For a set $X$ we define the function
\begin{center}
\xext=1000 \yext=100
\begin{picture}(\xext,\yext)(\xoff,\yoff)
\putmorphism(0,20)(1,0)[\sum_{n\in\o}\left[\begin{array}{c} X \\ n  \end{array}\right] \otimes_n A_n`F(X) `\tau_X]{1200}{1}a
\end{picture}
\end{center}
by
\[ \tau_X([\vec{x},a])=F(\vec{x})(a) \]
where $n\in\o$, $\vec{x}:(n]\ra X$ is an injection and $a\in A(n]$.

First we show that the transformation $\tau :\hat{A}\ra F$ is natural. Let $f:X\ra Y$ be a function. Then we have an epi/mono factorization of $f'$, $\vec{y}$ of the function $f\circ\vec{x}$
\begin{center} \xext=600 \yext=450
\begin{picture}(\xext,\yext)(\xoff,\yoff)
\setsqparms[2`3`1`2;600`400]
 \putsquare(0,0)[(n]`X`(m]`Y;\vec{x}`f'`f`\vec{y}]
 \end{picture}
\end{center}
Thus we have
\[ F(f)\tau_X([\vec{x},a]) =F(f\circ \vec{x})(a) = \]
\[ F(\vec{y}\circ f')(a) =F(\vec{y})(F(f')(a)) = \]
\[ =\tau_Y(A(f')(a)) \]
i.e. $\tau$ is natural.

It remains to show that $\tau$ is an isomorphism. Fix a set $X$. Let $x\in F(x)$. As $F$ is finitary, there is an $n\in\o$, $f:(n]\ra X$ and $y\in F(n]$ such that $F(f)(y)=x$. Let $\alpha$, $g$ be an epi/mono factorization of $f$
\begin{center} \xext=800 \yext=450
\begin{picture}(\xext,\yext)(\xoff,\yoff)
  \settriparms[1`3`-2;300]
  \putVtriangle(0,0)[(n]`X`{(m]};f`\alpha`g]
\end{picture}
\end{center}
Then $[g,A(\alpha)(y)] \in \left[\begin{array}{c} X \\ m  \end{array}\right] \otimes_m A_m$ and
\[ \tau_X([g,A(\alpha)(y)]) = F(g)(A(\alpha)(y)) = \]
\[ = F(g\circ \alpha)(y) = F(f)(y)= x\]
Thus $\tau_X$ is onto.

Now suppose that $[\vec{x},a]\in \left[\begin{array}{c} X \\ n  \end{array}\right] \otimes_n A_n$ and  $[\vec{y},b]\in \left[\begin{array}{c} X \\ m  \end{array}\right] \otimes_m A_m$ and
\[ \tau_X([\vec{x},a])  = \tau_X([\vec{y},b]) \]
We form a pullback of monos
\begin{center} \xext=600 \yext=500
\begin{picture}(\xext,\yext)(\xoff,\yoff)
\setsqparms[2`-2`-2`2;600`400]
 \putsquare(0,50)[(n]`X`(k]`{(m]};\vec{x}`f`\vec{y}`g]
 \end{picture}
\end{center}
that $F$ preserves, i.e. we have a pullback
\begin{center} \xext=1200 \yext=1050
\begin{picture}(\xext,\yext)(\xoff,\yoff)
\setsqparms[2`-2`-2`2;600`400]
 \putsquare(550,550)[F(n]`F(X)`F(k]`{F(m]};F(\vec{x})`F(f)`F(\vec{y})`F(g)]

\putmorphism(0,950)(1,0)[A(n]`\phantom{F(n]}`]{500}{2}a
\putmorphism(1150,550)(0,-1)[\phantom{F(m]}`A(m]`]{500}{-2}r
 \end{picture}
\end{center}
Since
\[ F(\vec{x})(a)= \tau_X([\vec{x},a]) = \tau_X([\vec{y},b]) = F(\vec{y})(b) \]
there is a $c\in F(k]$ such that
\[ F(f)(c)=a, \hskip 5mm F(g)(c)=b \]
As $a\in A(n]$ and $b\in B(m]$ we must have that $n=k=m$ and both $f$ and $g$ are  bijections. Put $\sigma=f\circ g^{-1}$. Then
\[ \vec{x}\circ\sigma = \vec{y}, \hskip 5mm a=A(\sigma)(b) \]
which means that $[\vec{x},a] = [\vec{y},b]$. This $\tau_X$ is an injection, as required. $\boxempty$
\vskip 2mm

A monad  $(M,\eta,\mu)$ on $Set$ is a {\em semi-analytic monad} iff $M$ is a semi-analytic functor and both $\eta$ and $\mu$
are semi-cartesian natural transformations. The category of semi-analytic monads with the semi-cartesian morphisms will be denoted by  $\SanMnd$.

We have inclusion functors
\begin{center}
\xext=1800 \yext=100
\begin{picture}(\xext,\yext)(\xoff,\yoff)
\putmorphism(0,20)(1,0)[\AMnd`\SanMnd`\cP_r^m]{800}{1}a
\putmorphism(800,20)(1,0)[\phantom{\SanMnd}`\Mnd`\cP_f^m]{800}{1}a
\end{picture}
\end{center}

%\newpage
\section{Full operads}\label{full_operads}
In this section we describe the equivalence of the category of full operads $\FOp$ to categories of equational theories $\ET$, Lawvere theories $\LT$ and monads $\Mnd$. As we said, the categories $\FOp$ and $\Mnd$ are categories of monoids of two equivalent monoidal categories thus they are obviously equivalent. Therefore we will just define the equivalences of categories
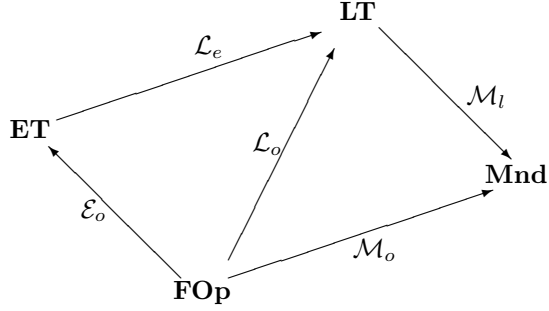
\begin{center} \xext=1800 \yext=1250
\begin{picture}(\xext,\yext)(\xoff,\yoff)
%level f
 %cats
  \put(620,40){$\FOp$}
  \put(1800,480){$\Mnd$}
  \put(0,650){$\ET$}
  \put(1250,1100){$\LT$}
%fun
    \put(920,600){$\cL_o$}
  \put(1300,200){$\cM_o$}
  \put(270,350){$\cE_o$}
  \put(700,970){$\cL_e$}
  \put(1730,780){$\cM_l$}

 %horizontal
 \put(830,120){\vector(3,1){1000}}
 \put(180,720){\vector(3,1){1000}}
 \put(650,120){\vector(-1,1){500}}
 \put(1400,1070){\vector(1,-1){500}}
 \put(830,190){\vector(1,2){400}}
\end{picture}
\end{center}
with the domain $\FOp$ and their essential inverses leaving all the verifications to the reader.  In \cite{Tr} it was shown that the category of full operads is equivalent to the category of abstract clones. The latter category in known to be equivalent to $\ET$.

Let us fix a morphism $h:\cO\ra \cO'$ of $\bF$-operads. We denote by $\kappa^n_i: (1]\ra (n]$ the function that sends $1$ to $i$, where $n\in \o$ and we denote  $\pi^n_i$  the action of $\kappa^n_i$ on the unit $\iota\in \cO_1$, i.e. $\pi^n_i=\kappa^n_i\cdot \iota$.

\subsection*{The functor $\cE_o$} The equational theory $\cE_o(\cO)=(L,A)$ has as the set of operation  $L=\coprod_{n\in\o} \cO_n$. For $f\in \cO_k$, $f_i\in \cO_{n_i}$, $i=1,\ldots, k$, $n=\sum_{i=1}^k n_i$, and $\psi: (k]\ra (n]$ we have the following axioms in $A$:
\[ f(f_1(x_1,\ldots, x_{n_1}),\ldots, f_k(x_{n-n_k+1}, \ldots, x_n))= ((f_1,\ldots,f_k)\star f)(x_1,\ldots, x_n) : \vec{x}^n \]

\[  f(x_{\psi(1)}, \ldots, x_{\psi(k)}) = (\psi \cdot f)(x_1, \ldots x_n) : \vec{x}^n \]
and
\[ \iota(x_1) = x_1 : \vec{x}^1 \]
The interpretation $\cE_o(h): \cE_o(\cO)\ra \cE_o(\cO')$ sends function symbol $f\in \cO_n$ to the term $h(f)(x_1,\ldots, x_n):\vec{x}^n$. This ends the definition of the functor $\cE_o$.

The essential inverse of the functor $\cE_o$
\[ \Op_e : \ET\lra \FOp \]
is defined as follows. For an equational theory $T=(L,A)$, the $\F$-operad $\cQ=\Op_e(T)$ has as $n$-ary operations $\cQ_n$ the terms in context $\vec{x}^n$
modulo provability in $T$. The action
\[ \cdot : \F(n,m) \times \cQ_n \lra \cQ_m \]
is defined as
\[ \phi\cdot [t(x_1, \ldots, x_n):\vec{x}^n] = [ t(x_{\phi(1)}, \ldots, x_{\phi(n)}):\vec{x}^m] \]
where $[t(x_1, \ldots, x_n):\vec{x}^n]$ is an operation in a $\cQ_n$ and $\phi : (n]\ra (m]$ is a function. The composition in $\cQ$
\[ \star : \cQ_{n_1}\times \ldots\times \cQ_{n_k}\times \cQ_k \lra Q_n \]
where $n=\sum_{i=1}^k n_i$ is defined by substitution with $\alpha$-conversion, i.e.
for $[t_i(x_1, \ldots, x_{n_i}):\vec{x}^{n_i}]\in \cQ_{n_i}$ and $[s(x_1, \ldots, x_k):\vec{x}^k]\in \cQ_k$ we have
\[ ([t_1(x_1, \ldots, x_{n_1}):\vec{x}^{n_1}],\ldots ,[t_1(x_1, \ldots, x_{n_k}):\vec{x}^{n_k}])\star [s(x_1, \ldots, x_k):\vec{x}^k] = \]
\[ =[s(t_1(x_{\alpha_1(1)}),\ldots,x_{\alpha_1(n_1)})), \ldots,   t_k(x_{\alpha_k(1)}),\ldots,x_{\alpha_k(n_k)})):\vec{x}^k] \]
where $\alpha_i : (n_i]\ra (n]$, for $i=1,\ldots, k$ is the obvious embedding.

\subsection*{The functor $\cL_o$} A morphism
\[ \lk f_1,\ldots, f_k\rk :n\lra k \]
 in the Lawvere theory $\cL_o(\cO)$ is a $k$-tuple such that $f_i\in \cO$, for $i=1,\ldots, k$. The identity on $n$ is
\[ \lk \pi^n_1,\ldots, \pi^n_n\rk :n\lra n \]
and the $i$-th projection from $n$ is
\[ \pi^n_i :n\lra 1 \]
Recall that $\pi^n_i$ is the value of the action of the function $\kappa^n_i: (1]\ra (n]$ that picks $i$ on the unit $\iota$.
The composition of morphisms
\[  \lk g_1,\ldots, g_n\rk :m\lra n, \;\;\;\;  \lk f_1,\ldots, f_k\rk :n\lra k \]
is
\[ \lk t_1,\ldots, t_k\rk :m\lra k \]
where
\[ t_i = \psi \cdot ((g_1, \ldots g_n)\ast f_i) \]
for $i=1,\ldots, k$, where $\psi : (n\cdot m] \ra (m]$ is given by $$\psi(r)= ((r-1)\, {\rm mod}\, m)+1$$ for $r\in  (n\cdot m]$.
The interpretation $\cL_o(h): \cL_o(\cO)\ra \cL_o(\cO')$ sends morphism  $\lk f_1,\ldots, f_k\rk :n\lra k$ to morphism $\lk h(f_1),\ldots, h(f_k)\rk :n\lra k$.
This ends the definition of the functor $\cL_o$.

The essential inverse of the functor $\cL_o$
\[ \Op_l : \LT\lra \FOp \]
is defined as follows. For a Lawvere theory $\bT$ the set of $n$-ary operation of the operad $Q=\Op_l(\bT)$ is $Q_n=\bT(n,1)$. The action in $\cQ$ is given by
\[ \phi \cdot f = f\circ \pi_\phi \]
where $f\in \cQ_n$, $\phi$ as above, and $\pi:\F^{op}\ra \bT$ is the unique morphism from the initial Lawvere theory. The composition in $\cQ$ is given by products and composition in $\bT$
\[ (f_1,\ldots, f_k)\star g = g\circ (f\times\ldots, \times f_k) \]

\subsection*{The functor $\cM_o$}  Let $f:X\ra Y$ be a function. The functor part of the monad $\cM_o(\cO)$ is defined via the coend formula
\[ \cM_o(\cO)(X) = \int^{n\in \F} X^n \times \cO_n\]
Thus an element of this set is an equivalence class of pairs $\lk \vec{x},t\rk \in X^n\times O_n$. We put
\[ \cM_o(\cO)(f)([\vec{x},t]) = [f\circ\vec{x},t] \]
The unit of the monad
\[ \eta^\cO_X : X \ra \cM(\cO)(X) \]
is given by
\[ x \mapsto [\tilde{x},\iota] \]
where $\tilde{x} : (1] \ra X$ sends $1$ to $x$. The multiplication
\[ \mu^\cO_X: \cM_o^2(\cO)(X) = \int^{k,n_1,\ldots,n_k\in \o} X^n\times \cO_{n_1}\times \ldots \times \cO_{n_k}\times \cO_k \lra  \cM_o(\cO)(X) \]
where $n=\sum_{i=1}^kn_i$, is given by
\[ [\vec{x},s_1,\ldots,s_k,t] \mapsto [\vec{x},(s_1,\ldots,s_k)\ast t] \]

The essential inverse of the functor $\cM_o$
\[ \Op_m : \Mnd\lra \FOp \]
is defined as follows. For a monad $T=(T,\eta,\mu)$ in $\Mnd$ the set of $n$-ary operations in the operad $\cQ=Op_m(T)$ is $Q_n=T(n]$. The action on $t\in Q_n$ is
\[ \phi\cdot t = T(\phi)(t) \]
Let $t_i\in T(n_i]$ and $s\in T(k]$. then we have corresponding morphisms
\[ \lceil s \rceil: (1]\lra T(k],\;\;\;  \lceil \vec{t} \rceil: (k]\lra T(n_1]+\ldots t(n_k] \]
'picking' $s$ and the $k$ terms $t_1,\ldots, t_k$. The action
\[ (t_1,\ldots,  t_k)\star s \]
is the element 'picked' by the morphism
\begin{center} \xext=3600 \yext=250
\begin{picture}(\xext,\yext)(\xoff,\yoff)
  \putmorphism(0,50)(1,0)[(1]`T(k]`\lceil s \rceil]{400}{1}a
  \putmorphism(400,50)(1,0)[\phantom{T(k]}`T(T(n_1]+\ldots + T(n_k])`T(\lceil\vec{t}\rceil)]{1000}{1}a
  \putmorphism(1400,50)(1,0)[\phantom{T(T(n_1]+\ldots + T(n_k])}`{T^2(n]}`{T({[T(\alpha_1), \ldots, T(\alpha_k)]})}]{1700}{1}a
  \putmorphism(3100,50)(1,0)[\phantom{{T^2(n]}}`{T(n]}`\mu]{500}{1}a
\end{picture}
\end{center}
in $T(n]$, where functions $\alpha_i$ are as above.

%\newpage
\section{Regular Lawvere theories vs regular operads}\label{sec_Lawvere}

In this section we study the relations between Lawvere theories and regular operads.  We shall describe the adjunction  $\Q^{ol}_r\dashv \P^{lo}_f$ and the properties of the embedding $\P^{lo}_f$.
\begin{center} \xext=1100 \yext=900
\begin{picture}(\xext,\yext)(\xoff,\yoff)
\setsqparms[1`-1`1`1;1100`800]
 \putsquare(0,50)[\FOp`\LT`\RegOp`\RegLT;\cL^f_o`\P_f^o`\Q^l_r`\cL^r_o]
 \put(50,780){\vector(0,-1){650}}
 \put(70,430){$\Q^{o}_r$}

   \put(150,150){\vector(4,3){800}}
   \put(420,500){$\P^{ol}_f$}

   \put(990,730){\vector(-4,-3){800}}
   \put(580,360){$\Q^{lo}_r$}

 \put(1050,130){\vector(0,1){650}}
 \put(920,430){$\P^{l}_f$}
\end{picture}
\end{center}

\subsection*{The functor $\P_f=\P_f^{ol} : \RegOp\ra \LT$}
Let $\cO$ be a regular operad; $\iota$, $\cdot$, $\ast$ denote the unit, action (of $\S$), and composition in $\cO$, respectively. We define a Lawvere theory $\P_f(\cO)$ as follows. The set of objects of $\P_f(\cO)$ is the set of natural numbers $\o$. A morphism from $n$ to $m$ is an equivalence class of spans
\begin{center} \xext=800 \yext=520
\begin{picture}(\xext,\yext)(\xoff,\yoff)
 \settriparms[1`1`0;400]
  \putAtriangle(0,0)[r`n`m;\phi`\lk f,g_i\rk_{i\in m}`]
\end{picture}
\end{center}
such that $\phi:(r]\ra(n]$ is a function called {\em amalgamation}, $f:(r]\ra(m]$ is a monotone function called the {\em arity function } (as it determines the arities of the operations $g_i$), $r_i=|f^{-1}(i)|$ and we have $g_i\in \cO_{r_i}$ for $i\in (m]$ and $r=\sum_{i=1}^m r_i$, moreover $\phi$ and $f$ are jointly mono (or what comes out the same, $\phi$ is mono on the fibers of $f$). Two spans $\lk \phi,f,g_i\rk_{i\in m}$ and $\lk \phi ',f',g'_j\rk_{j\in m'}$ are equivalent
whenever $f=f'$ and there are permutations $\sigma_i : (r_i]\ra (r_i]$ for $i\in(m]$
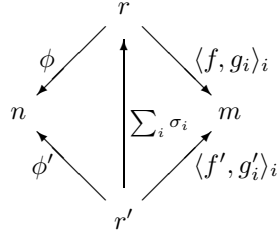
\begin{center} \xext=800 \yext=920
\begin{picture}(\xext,\yext)(\xoff,\yoff)
 \settriparms[1`1`0;400]
  \putAtriangle(0,400)[r`n`m;\phi`\lk f,g_i\rk_{i}`]
  \settriparms[0`-1`-1;400]
  \putVtriangle(0,0)[\phantom{n}`\phantom{m}`r';`\phi'`{\lk f',g'_i\rk_{i}}]
  \putmorphism(400,750)(0,-1)[``]{700}{-1}a
    \put(420,380){$_{\sum_i \sigma_i}$}
\end{picture}
\end{center}
such that
\[  g_i  = \sigma_i \cdot g'_i,\;\;\;\; \phi\circ \sum_i \sigma_i =\phi'\]
In particular $|f^{-1}(i)|=r_i=|f'^{-1}(i)|$, for $i\in (m]$.
By $\sum_i \sigma_i:r \ra r$ we mean the permutation that is formed by placing permutations $\sigma_i$ 'one after another'. Thus, it respects the fibers of $f$, i.e.  $f\circ \sum_i \sigma_i =f'$. A morphism in the category $\P_f(\cO)$ is a class of spans modulo the above equivalence relation.

We could represent morphisms in $\P_f(\cO)$ as spans without the requirement that $\phi$ and $f$ are jointly mono. But then the relation identifying the spans would be more complicated. Instead of permutations $\sigma_i$ we would need to consider surjections. But this relation is not an equivalence relation and we would need to work with the equivalence relation generated by such a relation. However, it might happen, as with the compositions defined below, that an operation on spans naturally leads to a span whose amalgamation $\phi$ is not injective on fibers of $f$. In such a case we can {\em regularize} the span as follows. Let  $\lk \phi,f,g_i\rk_{i\in m}$ be a span as above but with $\phi$ not necessarily injective on the fibers of $f$. Let $\phi_i$ be the restriction of $\phi$ to the fiber $f^{-1} (i)$ for $i\in (m]$. Let $\phi_i=\phi_i'\circ \psi_i$ be an epi-mono factorization of $\phi_i$  and $g_i'= \psi_i\cdot g_i$, for $i\in (m]$. Then the {\em regularization} of the span  $\lk \phi,f,g_i\rk_{i\in m}$ is the span
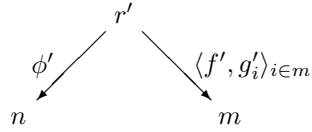
\begin{center} \xext=800 \yext=520
\begin{picture}(\xext,\yext)(\xoff,\yoff)
 \settriparms[1`1`0;400]
  \putAtriangle(0,0)[r'`n`m;\phi'`\lk f',g'_i\rk_{i\in m}`]
\end{picture}
\end{center}
where $r'=\sum_{i\in (m]} r_i'$, $r_i'=|dom(\phi_i')|$, for $i\in (m]$,  $\phi'=(\phi_1'+\ldots+\phi_m'):r'\ra n$. $f'$ is the monotone map sending the elements in the domain of $ \phi_i'$ to $i$, for $i\in (m]$.

The composition of morphisms $\lk \phi,f,g_i\rk_{i\in (m]}: n \ra m$ and $\lk \phi ',f',g'_j\rk_{j\in (k]}:m\ra k$ to $\lk \phi'',f'',g''_j\rk_{i\in (k]}: n \ra k$ is
defined in two steps as follows. In the diagram
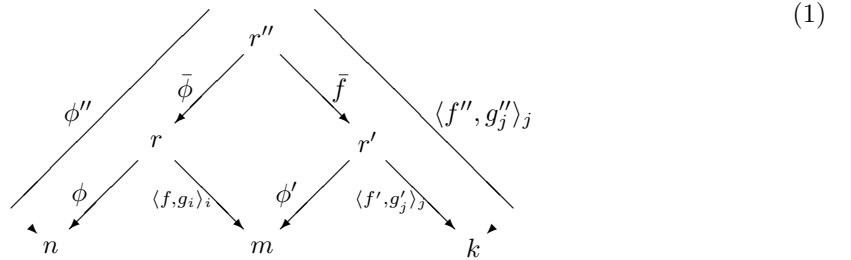
\begin{equation}\label{compo}\end{equation}
\begin{center} \xext=2000 \yext=720
\begin{picture}(\xext,\yext)(\xoff,\yoff)
 \settriparms[1`1`0;400]
  \putAtriangle(200,0)[r`n`m;\phi``]
   \putAtriangle(1000,0)[r'`\phantom{m}`k;\phi'``]
   \put(580,180){$_{\lk f,g_i\rk_{i}}$}
    \putAtriangle(600,400)[r''`\phantom{r}`\phantom{r'};\bar{\phi}`\bar{f}`]
   \put(1350,180){$_{\lk f',g'_j\rk_{j}}$}

   \put(800,900){\line(-1,-1){750}}
   \put(800,900){\line(1,-1){100}}
   \put(50,150){\vector(1,-1){100}}
   \put(250,480){$\phi''$}

    \put(1200,900){\line(1,-1){750}}
   \put(1200,900){\line(-1,-1){100}}
   \put(1950,150){\vector(-1,-1){100}}
   \put(1650,480){$\lk f'',g''_j \rk_j$}
\end{picture}
\end{center}
the square is a pullback of $f$ along $\phi'$. The function $\bar{f}$ is so chosen that it is monotone. We define a span by $f''=f'\circ \bar{f}$,  $\phi''=\phi'\circ \bar{\phi}$, and $g_j''= g_j'\ast \lk g_{\phi'(l)} \rk_{l\in f'^{-1}(j)}$. Finally, we take a regularization of this span to get a regular span that represents the composition. We leave for the reader the verification that the composition is a congruence with respect to the equivalence relation on spans.

The identity on $n$ is the span
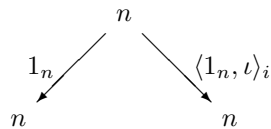
\begin{center} \xext=800 \yext=500
\begin{picture}(\xext,\yext)(\xoff,\yoff)
 \settriparms[1`1`0;400]
  \putAtriangle(0,0)[n`n`n;1_n`\lk 1_n,\iota\rk_{i}`]
\end{picture}
\end{center}

The $i$-projection $\pi^n_i:n\ra 1$ on $i$-th coordinate is the span
\begin{center} \xext=800 \yext=520
\begin{picture}(\xext,\yext)(\xoff,\yoff)
 \settriparms[1`1`0;400]
  \putAtriangle(0,0)[1`n`1;\tilde{i}`\lk 1_1,\iota\rk`]
\end{picture}
\end{center}
where $i\in (n]$ and $\tilde{i}(1)=i$.

For a morphism of regular operads $h:\cO\ra \cO'$, we define a functor
\[ \P_f(h): \P_f(\cO) \lra \P_f(\cO') \]
so that for a morphism $\lk \phi,f,g_i\rk_{i\in (m]}: n \ra m$ in $\P_f(\cO)$, we have
\[ \P_f(h)(\lk \phi,f,g_i\rk_{i\in (m]}) = \lk \phi,f,h(g_i)\rk_{i\in (m]}: n \ra m \]
in  $\P_f(\cO')$. This ends the definition of the functor $\P_f$.
\subsection*{The functor $\Q_r=\Q_r^{lo} : \LT \lra \RegOp$}

Let $\bT$ be a Lawvere theory. Recall that $\pi:\F^{op}\ra \bT$ is the morphism from the initial Lawvere theory. The operad $\Q_r(\bT)$ consists of operations of $\bT$, i.e. morphisms to $1$. In detail it can be described as follows. The set of $n$-ary operations $\Q_r(\bT)_n$ is the set of  $n$-ary operations  $\bT(n,1)$ of $\bT$, for $n\in \o$. The action
\[ \cdot : \S(n,m) \times \Q_r(\bT)_n \lra \Q_r(\bT)_m \]
is given, for $f\in\bT(n,1)$ and $\phi\in \S(n,m)$,  by
\[ \phi \cdot  f = f\circ \pi_{\phi} \]
The identity of $\Q_r(\bT)$ is $\iota=1_1\in\bT(1,1)$.
The composition
\[ \ast : \Q_r(\bT)_{n_1}\times\Q_r(\bT)_{n_{k}}\times\Q_r(\bT)_k\lra \Q_r(\bT)_n \]
is given, for $f\in  \Q_r(\bT)_k$ and $f_i\in \Q_r(\bT)_{n_i}$,  where $i\in (k]$, $n={\sum_{i\in k} n_i}$, by
\[  \lk f_1,\ldots , f_{k}\rk \ast f = f\circ (f_1\times\ldots ,\times f_{k}) \]
where $f_1\times\ldots ,\times f_{k}$ is defined using the chosen projections in $\bT$ and $\circ$ is the composition in $\bT$.

If $F:\bT\ra \bT'$ is a morphism of Lawvere theories then the map of regular operads
\[ \Q_r(F) : \Q_r(\bT) \ra \Q_r(\bT') \]
is defined, for $f\in \Q_r(\bT)_n$, by
\[ \Q_r(F)(f) = F(f) \]
This ends the definition of the functor $\Q_r$.

\subsection*{The adjunction $\P_f\dashv \Q_r$ and the properties of the functor $\P_f$}

We note for the record
\begin{proposition} \label{funtors_P_Op}
The functors $\P_f : \RegOp\lra \LT$ and $\Q_r : \LT \ra \RegOp$ are well defined. $\boxempty$
\end{proposition}

We have an easy

\begin{lemma} \label{iso}
Let $\cO$ be a regular operad and $n\in\o$. An isomorphism on $n$ in  $P_r(\cO)$ has a representation
by a span of the following form
\begin{center} \xext=800 \yext=500
\begin{picture}(\xext,\yext)(\xoff,\yoff)
 \settriparms[1`1`0;400]
  \putAtriangle(0,0)[n`n`n;\phi`\lk 1_n,a_i\rk_{i}`]
\end{picture}
\end{center}
where $\phi:(n]\ra (n]$ is a bijection, $a_i\in \cO_1$ is an invertible operation, i.e. there is $b_i\in \cO_1$ such that $a_i\ast b_i=\iota =b_i\ast a_i$ for $i\in (n]$. It is the unique span in its equivalence class.
\end{lemma}

{\em Proof.} Suppose that we have two spans
\begin{equation}\end{equation}
\begin{center} \xext=2000 \yext=720
\begin{picture}(\xext,\yext)(\xoff,\yoff)
 \settriparms[1`1`0;400]
  \putAtriangle(200,0)[r`n`n;\phi``]
   \putAtriangle(1000,0)[r'`\phantom{n}`n;\phi'``]
   \put(580,180){$_{\lk f,g_i\rk_{i}}$}
   \put(1750,180){$_{\lk f',h_j\rk_{j}}$}
\end{picture}
\end{center}
representing two  morphisms that compose to identity $1_n$ both ways. Then, since the displayed composition is $1_n$, the functions  $\phi$, $f'$ are surjections. Hence the functions  $\phi'$, $f$ are surjections, as well. Having this it is easy to notice that $\phi$ sends elements in different fibers of $f$ to different elements. Thus $\phi$ (and $\phi'$)  must be in fact a bijection.

If we take a composition $\lk g_i\rk_{i\in f'^{-1}(j)}\ast h_j$ and regularize it multiplying by a surjection, we will get $\iota$. This implies that
$h_j$ cannot be a nullary operation and that the arity of $g_i$'s is at most one. Now, as the above spans represent morphisms that compose both ways to identity, it is easy to see that we get the required description.
 $\boxempty$

\begin{proposition} \label{operad-emb}
We have an adjunction $\P_f\dashv \Q_r$. The functor $\P_f$ is faithful.
\end{proposition}

{\em Proof.}
 We shall show that $\P_f\dashv \Q_r$. For a regular operad $\cO$ the unit is
\[ \eta_\cO :\cO \lra \Q_r (\P_f(\cO)) \]
\[ \cO_n \ni g \mapsto \lk 1_n,!,g\rk  \]
For Lawvere theory $\bT$ the counit is
\[ \varepsilon_\bT: \P_f\Q_r(\bT)\lra \bT   \]
\[ \lk \phi,f, g_i\rk_{i\in (m]} \mapsto (g_1\times \ldots \times g_{m})\circ \pi_\phi  \]

We verify the triangular equalities. For $g\in \Q_r(\bT)_n=\bT(n,1)$ we have
\[ \Q_r (\varepsilon_\bT) \circ \eta_{\Q_r(\bT)}(g) = \]
\[  =\Q_r (\varepsilon_\bT) (\lk 1_n,!,g\rk) = \]
\[ = g\circ \pi_{1_n} = g\]
For $\lk \phi,f,g_i\rk_{i\in (m]}\in \P_f(\cO)$ we have
\[ \varepsilon_{\P_f(\cO)} \circ \P_f(\eta_\cO)(\lk \phi,f,g_i \rk_{i\in (m]} ) = \]
\[ =\varepsilon_{\P_f(\cO)} (\lk \phi,f,\lk 1_{r_i},!,g_i\rk \rk_{i\in (m]} ) = \]
\[ = (\lk 1_{r_1},!,g_1 \rk\times \ldots \times\lk 1_{r_{m}},!,g_{m} \rk) \circ \pi_\phi  =\]
\[ = \lk \phi,f,g_i \rk_{i\in (m]}  \]
As the unit $\eta_\cO$ is mono, $\P_f$ is faithful.
$\boxempty$

\begin{corollary} \label{operad-emb-im}
The triangle
\begin{center} \xext=1100 \yext=900
\begin{picture}(\xext,\yext)(\xoff,\yoff)
\setsqparms[1`-1`0`0;1100`800]
 \putsquare(0,50)[\FOp`\LT`\RegOp`;\cL^f_o`\P_f^o``]

   \put(150,150){\vector(4,3){800}}
   \put(420,500){$\P^{ol}_f$}
\end{picture}
\end{center}
commutes up to an isomorphism.
\end{corollary}

{\em Proof.} By Proposition \ref{operad-emb} it is enough to show that triangle
\begin{center} \xext=1100 \yext=1000
\begin{picture}(\xext,\yext)(\xoff,\yoff)
\setsqparms[-1`0`0`0;1100`800]
 \putsquare(0,50)[\FOp`\LT`\RegOp`;\Op^f_l```]
 \put(50,780){\vector(0,-1){650}}
 \put(70,430){$\Q^{o}_r$}

   \put(990,730){\vector(-4,-3){800}}
   \put(580,360){$\Q^{lo}_r$}
\end{picture}
\end{center}
commutes up to an isomorphism, where the functor $\Op^f_l$ is the left adjoint to $\cL^f_o$ and together they form an adjoint equivalence. The functor $\Op^f_l$ is defined as the functor $\Q_r$ except the action involved is the action of the whole category $\F$ rather than its subcategory $\S$. As the functor $\Q^{o}_r$ forgets this additional part of the action the above diagram clearly commutes.
$\boxempty$

\begin{proposition}
The functor $\P_f: \RegOp\ra \LT$ is  full on isomorphisms and its essential image is $\RegLT$. In particular $\RegOp$ is equivalent to $\RegLT$.
\end{proposition}

{\em Proof.} Recall that we have a unique morphism of Lawvere theories from the initial theory $\pi : \mathds{F}^{op}\ra\P_f(\cO)$.
For a function $\phi : m\ra n$, $\pi_\phi$ the morphism $\pi_\phi$ is represented by the span of form
\begin{center} \xext=800 \yext=520
\begin{picture}(\xext,\yext)(\xoff,\yoff)
 \settriparms[1`1`0;400]
  \putAtriangle(0,0)[m`n`m;\phi`\lk 1_m,\iota\rk_{i\in (m]}`]
\end{picture}
\end{center}

The class of projection morphisms, is the closure under isomorphism of the class of morphisms $\{\pi_\phi : \phi\in \I  \}$ in $\P_f(\cO)$. Using Lemma \ref{iso}, it is easy to see that the projection morphisms in $\P_f(\cO)$ are (represented by) the spans of the form
\begin{center} \xext=800 \yext=520
\begin{picture}(\xext,\yext)(\xoff,\yoff)
 \settriparms[1`1`0;400]
  \putAtriangle(0,0)[m`n`m;\phi`\lk 1_m,a_i\rk_{i\in (m]}`]
\end{picture}
\end{center}
where $\phi$ is an injection and $a_i$ is an invertible unary operation.

The regular morphisms in $\P_f(\cO)$ are (represented by) the
spans of the form
\begin{center} \xext=800 \yext=520
\begin{picture}(\xext,\yext)(\xoff,\yoff)
 \settriparms[1`1`0;400]
  \putAtriangle(0,0)[n`n`m;\phi`\lk f,g_i\rk_{i\in (m]}`]
\end{picture}
\end{center}
where $\phi$ is a bijection.

Clearly, both classes contain isomorphisms and are closed under composition.

Any morphism $\lk \phi,f,g_i\rk_{i\in (m]}:n\ra m$ in $\P_f(\cO)$ has a projection-regular factorization as follows
\begin{center} \xext=2000 \yext=520
\begin{picture}(\xext,\yext)(\xoff,\yoff)
 \settriparms[1`1`0;500]
  \putAtriangle(0,0)[r`n`r;\phi'``]
  \putAtriangle(1000,0)[r`\phantom{r}`n;1_r`\lk f,g_i\rk_{i}`]
  \put(430,200){$\lk 1_{r},\iota\rk_{j}$}
\end{picture}
\end{center}

Thus to show that projections and (what we have described as) regular morphisms form a factorization system it remains to show that projection morphisms are left orthogonal to the regular morphisms. Let
\begin{center} \xext=1250 \yext=1250
\begin{picture}(\xext,\yext)(\xoff,\yoff)
%lt
\setsqparms[-1`-1`0`0;600`500]
 \putsquare(50,600)[n`r`\phantom{k}`;\psi`\phi``]
%rt
\setsqparms[1`0`-1`0;600`500]
 \putsquare(650,600)[\phantom{r}`m``m;\lk f,h_i \rk_{i\in (m]}``1_m`]
%lb
\setsqparms[0`1`0`-1;600`500]
 \putsquare(50,100)[k``k`r';`\lk 1_k,a_j \rk_{j\in (k]}``\phi']
%rb
\setsqparms[0`0`1`1;600`500]
\putsquare(650,100)[`\phantom{m}`\phantom{r'}`1;``\lk !,g \rk`\lk !,g' \rk]

\putmorphism(650,1050)(0,-1)[``\sigma]{900}{1}r
\end{picture}
\end{center}
be a commutative square in $\P_f(\cO)$ with the left vertical morphism $\lk \phi,1_r,a_i\rk_{j\in (k]}$ a structural map and right vertical morphism  $\lk 1_m,!,g\rk$ a regular map. We have chosen the right bottom to be $1$ to simplify notation but the general case is similar. The commutativity means that $r=r'$ and there is a permutation $\sigma\in S_r$  such that
\[ \psi=\phi\circ \phi'\circ \sigma\]
and
\[\lk a_{\phi'(1)},\ldots,a_{\phi'(r)} \rk\ast g'\ = \sigma\cdot (\lk h_1,\ldots,h_{m} \rk\ast  g)\]
Putting into the square a diagonal morphism $\lk \phi'\circ \sigma,f,\bar{h_i}\rk_{i\in (m]}$
\begin{center} \xext=1250 \yext=1250
\begin{picture}(\xext,\yext)(\xoff,\yoff)
%lt
\setsqparms[-1`-1`1`0;600`500]
 \putsquare(50,600)[n`r`\phantom{k}`r;\psi`\phi`1_r`]
%rt
\setsqparms[1`0`-1`0;600`500]
 \putsquare(650,600)[\phantom{r}`m``m;\lk f,h_i \rk_{i}``1_m`]
%lb
\setsqparms[0`1`1`-1;600`500]
 \putsquare(50,100)[k`\phantom{r}`k`r;`\lk 1_k,a_j \rk_{j}`\sigma`\phi']
%rb
\setsqparms[0`0`1`1;600`500]
 \putsquare(650,100)[`\phantom{m}`\phantom{r}`1;``\lk !,g \rk`\lk !,g' \rk]

  \put(200,400){$_{\phi'\circ \sigma}$}
  \put(850,720){$_{\lk f,\bar{h_i} \rk_{i}}$}
  \put(600,550){\vector(-4,-3){500}}
   \put(700,650){\vector(4,3){500}}
\end{picture}
\end{center}
where
\[ \bar{h_i} = \lk a^{-1}_{\phi'\circ \sigma(l)} \rk_{l\in f^{-1}(i)}\ast h_i  \]
we see that the permutations $1_r$ and $\sigma$ show that both triangles commute. Thus regular morphisms are indeed right orthogonal to the structural ones and
$\P_f(\cO)$ is a regular Lawvere theory.

From the description of the functor  $\P_f(h): \P_f(\cO)\ra \P_f(\cO')$ and the description of the structure of $\P_f(\cO)$, it is clear that $\P_f(h)$ sends the regular (projection) morphisms to the regular (projection) morphisms. Thus $\P_f(h)$ is a regular morphism of Lawvere theories.

Now let $\bT$ be any Lawvere theory. As the class of regular morphisms in $\bT$ is right orthogonal to a class of morphisms and it is closed under composition, finite products and isomorphisms. Moreover, for any surjection in $\psi:(n]\ra (m]$ in $\F$ the image $\pi_\psi:m\ra n$ in $\bT$ is regular as it is orthogonal to all projection morphisms in $\bT$. Thus surjections in $\S$ act on all regular morphisms $f:n\ra 1$ in $\bT$ on the left
\[ \cdot: \S(n,m) \times \bT(n,1)\lra \bT(m,1) \]
by
\[ (\psi,f) \mapsto f\circ \pi_\psi \]
Hence the regular operations of any Lawvere theory $\bT$ form a regular operad. The unit is the  identity morphism on $1$. The composition $\lk f_1, \ldots, f_n\rk \ast f$ is defined to be $f\circ (f_1\times \ldots \times f_n)$. The action of $\S$ is defined as above. Let us denote this operad by $\bT^r$. We have an inclusion morphism of regular operads $\bT^r\ra \Q_r(\bT)$. By adjunction we get a morphism of Lawvere theories
\[ \xi_\bT:\P_f(\bT^r) \lra \bT \]
Clearly, $\xi_\bT$ is bijective on objects. If $\bT$ is regular then $\xi_\bT$ is full (faithful) since the projection-regular factorization exists (is unique and $\pi:\F^{op}\ra \bT$ is faithful).

If $I:\bT\ra \bT'$ is a regular interpretation between any Lawvere theories, then the diagram
\begin{center} \xext=1000 \yext=600
\begin{picture}(\xext,\yext)(\xoff,\yoff)
\setsqparms[1`1`1`1;1000`500]
 \putsquare(0,50)[\P_f(\bT^r)`\P_f(\bT'^r)`\bT`\bT;\P_f(I^r)`\xi_\bT`\xi_{\bT'}`I]
 \end{picture}
\end{center}
commutes, where $I^r$ is the obvious restriction of $I$ to $\bT^r$. Thus the essential image of $\P_f$ is indeed the category of regular Lawvere theories and regular interpretations. An isomorphic interpretation of Lawvere theories is always regular. Therefore $\P_f$ is full on isomorphisms. $\boxempty$

We have
\begin{proposition}\label{monadic}
The functor $\Q_r : \LT \ra \RegOp$ is monadic.
\end{proposition}
{\em Proof.}  We shall verify the assumption of the Beck monadicity theorem.
By Proposition \ref{operad-emb} $\Q_r$ has a left adjoint. It is easy to see that $\Q_r$
reflects isomorphisms. We shall verify that $\LT$ has and $\Q_r$ preserves $\Q_r$-contractible
 coequalizers.

Let $I,I': \bT' \ra \bT$ be a pair of interpretations between Lawvere theories so that
\begin{center} \xext=1400 \yext=400
\begin{picture}(\xext,\yext)(\xoff,\yoff)

\putmorphism(800,200)(1,0)[\phantom{\Q_r(\bT)}`\cO`q]{600}{1}a
\putmorphism(800,150)(1,0)[\phantom{\Q_r(\bT)}`\phantom{\cO}`s]{600}{-1}b
  \putmorphism(0,200)(1,0)[\Q_r(\bT')`\Q_r(\bT)`\Q_r(I')]{800}{1}b
   \putmorphism(0,50)(1,0)[\phantom{\Q_r(\bT')}`\phantom{\Q_r(\bT)}`r]{800}{-1}b
    \putmorphism(0,280)(1,0)[\phantom{\Q_r(\bT')}`\phantom{\Q_r(\bT)}`\Q_r(I)]{800}{1}a
\end{picture}
\end{center}
is a split coequalizer in $\RegOp$.  We define a Lawvere theory $\bT_\cO$ so that a morphism
from $n$ to $m$ in $\bT_\cO$ is an $m$-tuple $\lk g_1, \ldots, g_m\rk$ with $g_i\in \cO_n$, for $i=1,\ldots, m$.
The compositions and the identities in $\bT_\cO$ are defined in the obvious way from
the compositions and the unit in $\cO$. The projections $\bar{\pi}^n_i$ in $\bT_\cO$ are the images of the projections
$\pi^n_i$ in $\bT$, i.e. $\bar{\pi}^n_i= q(\pi^n_i)$.

The functor $\tilde{q}: \bT \ra \bT_\cO$ is defined, for $f:n\ra m$ in $\bT$, as
\[ \tilde{q}(f) = \lk q(\pi^m_1\circ f),\ldots, q(\pi^m_m\circ f)\rk \]

First we verify, that $\bT_\cO$ has finite products. For this, it is enough to verify that
$\lk f_1, \ldots, f_n\rk \ast\bar{\pi}^n_i = f_i$, where $\ast$ is the multiplication in the operad $\cO$.
The uniqueness of the morphism into the product is obvious from the construction. We have routine calculations
\[  \lk f_1, \ldots, f_n\rk \ast\bar{\pi}^n_i  = \]
\[  q\circ s( \lk f_1, \ldots, f_n\rk \ast q(\pi^n_i) ) = \]
\[  \lk q\circ s(f_1), \ldots, q \circ s(f_n)\rk\ast (q\circ s\circ q(\pi^n_i)) = \]
\[  \lk q\circ s(f_1), \ldots,q\circ s(f_n)\rk \ast (q\circ Op(I)\circ t(\pi^n_i)) = \]
\[  \lk q\circ s(f_1), \ldots,q\circ s(f_n)\rk  \ast (q\circ Op(I')\circ t(\pi^n_i))  = \]
\[   \lk q\circ s(f_1), \ldots,q\circ s(f_n)\rk \ast (q(\pi^n_i)) = \]
\[  q(\lk  s(f_1), \ldots,s(f_n)\rk\ast\pi^n_i) = \]
\[  q(s(f_i) = f_i \]

It is obvious that $\tilde{q}$ is a morphism of Lawvere theories and that $\Q_r(\tilde{q})=q$.
It remains to verify that $\tilde{q}$ is a coequalizer in $\LT$.
Let $p: \bT \ra \bS$ be a morphism in $\LT$ coequalizing $I$ and $I'$.
\begin{center} \xext=1100 \yext=700
\begin{picture}(\xext,\yext)(\xoff,\yoff)
 \settriparms[1`1`1;500]
  \putqtriangle(600,0)[\bT`\bT_\cO`\bS;\tilde{q}`p`\tilde{k}]

  \putmorphism(0,500)(1,0)[\bT'`\phantom{\bT}`]{600}{0}b
   \putmorphism(0,450)(1,0)[\phantom{\bT'}`\phantom{\bT}`I']{600}{1}b
    \putmorphism(0,550)(1,0)[\phantom{\bT'}`\phantom{\bT}`I]{600}{1}a
\end{picture}
\end{center}
The morphism $\Q_r(p)$ coequalizes  $\Q_r(I)$ and $\Q_r(I')$ in $\RegOp$.
Thus there is a unique morphism $k$ in $\RegOp$ making the triangle on the right
\begin{center} \xext=1300 \yext=700
\begin{picture}(\xext,\yext)(\xoff,\yoff)
 \settriparms[1`1`1;500]
  \putqtriangle(800,0)[\Q_r(\bT')`\cO`\Q_r(\bS);q`\Q_r(p)`k]

  \putmorphism(0,500)(1,0)[\Q_r(\bT')`\phantom{\Q_r(\bT)}`]{800}{0}b
   \putmorphism(0,450)(1,0)[\phantom{\Q_r(\bT')}`\phantom{\Q_r(\bT)}`\Q_r(I')]{800}{1}b
    \putmorphism(0,550)(1,0)[\phantom{\Q_r(\bT')}`\phantom{\Q_r(\bT)}`\Q_r(I)]{800}{1}a
\end{picture}
\end{center}
commute. We define the functor $\tilde{k}$ so that
\[  \tilde{k}(\lk f_1,\ldots, f_n\rk) = \lk k(f_1),\ldots, k(f_n)\rk \]
for any morphism $\lk f_1,\ldots, f_n\rk$ in $\bT_\cO$. The verification that  $\tilde{k}$ is the required unique functor
is left for the reader.
 $\boxempty$

%\newpage
\section{Semi-analytic monads vs  regular operads}\label{sec_mnd}

The main objective of this section is to show that the square
\begin{equation}\label{square_semi_an} \end{equation}
\begin{center} \xext=800 \yext=450
\begin{picture}(\xext,\yext)(\xoff,\yoff)
\setsqparms[1`-1`-1`1;800`500]
\putsquare(0,50)[\FOp`\Mnd`\RegOp`\SanMnd;\cM^f_o`\P^o_f`\P^m_f`\cM^r_o]
 \end{picture}
\end{center}
commutes up to isomorphism, with horizontal functors being equivalences of categories. The functor  $\cM^f_o$ was defined in Section \ref{full_operads}. $\P^m_f$ is an inclusion.

To this end we shall use Theorem \ref{im hat} and some general considerations but the functor part of the monad $\cM^r_o(\cO)$ will be described explicitly below. For a regular operad $\cO$, the monad $\cM^r_o(\cO)$ on a set $X$ is defined as follows. The functor $\cM^r_o(\cO)$ is the application of the functor $\hat {(-)} : Set^\S \lra End$ to $\cO$ considered as a functor $\cO:\S\ra Set$, see Section \ref{pres_of_cats_of_algs} for details.

In particular for a set $X$ we have
\[ \cM^r_o(\cO)(X)=\sum_{n\in \o} \left[\begin{array}{c} X \\ n  \end{array}\right]\otimes_{n} \cO_n \]
In the set $\left[\begin{array}{c} X \\ n  \end{array}\right]\otimes_{n} \cO_n$ we identify $\lk\vec{x}\circ \sigma,a\rk$ with $\lk \vec{x}, \sigma\cdot a\rk$ for $a\in \cO_n$, $\vec{x}: (n]\ra X$ and $\sigma\in S_n$.

Let $\gamma:\S\ra \F$ be the inclusion functor. It induces the following diagram of categories and functors that we describe below
\begin{center} \xext=3100 \yext=1300
\begin{picture}(\xext,\yext)(\xoff,\yoff)
%squares
\setsqparms[1`0`0`1;1600`750]
 \putsquare(150,450)[\hskip 25mm\Mnd\; =\; \Mon(\End)`\End`\hskip 25mm\SanMnd\; =\; \Mon(\San)`\San;\widehat{U}```U]
\setsqparms[-1`0`0`-1;800`750]
 \putsquare(1750,450)[\phantom{\End}`Set^\F`\phantom{\San}`Set^\S;i_\F```i_\S]

%vertical
 \putmorphism(700,1150)(0,-1)[\phantom{\End}`\phantom{\San}`\Mon((-)^{sa})]{700}{1}r
  \putmorphism(600,1150)(0,-1)[\phantom{\End}`\phantom{\San}`]{700}{-1}l

 \putmorphism(200,1150)(0,-1)[\phantom{\Mnd}`\phantom{2\bCat_{\times}}`\Q_r^m]{700}{1}r
 \putmorphism(100,1150)(0,-1)[\phantom{\Mnd}`\phantom{2\bCat_{\times}}`\P_f^m]{700}{-1}l

 \putmorphism(1800,1150)(0,-1)[\phantom{\End}`\phantom{\San}`(-)^{sa}]{700}{1}r
 \putmorphism(1700,1150)(0,-1)[\phantom{\End}`\phantom{\San}`i^{sa}]{700}{-1}l

  \putmorphism(2600,1150)(0,-1)[\phantom{Set^\F}`\phantom{Set^\S}`\gamma^*]{700}{1}r
 \putmorphism(2500,1150)(0,-1)[\phantom{Set^\F}`\phantom{Set^\S}`Lan_\gamma]{700}{-1}l

 \putmorphism(3100,1150)(0,-1)[\F`\S`\gamma]{700}{-1}r

%\Wb=\Mon(\W)
  \put(150,180){\oval(100,100)[b]}
  \put(100,180){\line(0,1){200}}
  \put(200,180){\vector(0,1){200}}
  \put(0,0){${(\Wb,\bar{\eta},\bar{\mu})=\Mon(\W,\eta,\mu) }$}
   \put(650,180){\oval(100,100)[b]}
  \put(600,180){\line(0,1){200}}
  \put(700,180){\vector(0,1){200}}

%\W
  \put(1500,0){${(\W,\eta,\mu)}$}
   \put(1750,180){\oval(100,100)[b]}
  \put(1700,180){\line(0,1){200}}
  \put(1800,180){\vector(0,1){200}}
 \end{picture}
\end{center}
$\gamma^*$ is the functor of composing with $\gamma$. It has a left adjoint $Lan_\gamma$, the left Kan extension along $\gamma$. For $C\in Set^\S$ it is given by the coend formula
\[ Lan_\gamma(C)(X) = \int^{n\in \S}  X^n\times C(n]     \]
The functor $i^{sa}: \San \ra \End$ is just an inclusion. The equivalence
\[  i_\S: Set^\S\lra \San  \]
is defined by left Kan extension that may be given by the coend formula and a coproduct
\[ i_\S(C)(X) \;\;=\;\; \int^{n\in \S} X^n\times C(n]\;\;=\;\;\sum_{n\in \o} \left[\begin{array}{c} X \\ n  \end{array}\right]\otimes_n C(n] \]
where $C\in Set^\S$. Similarly, the equivalence
\[  i_\F: Set^\F\lra \End  \]
is defined by left Kan extension that is given by the coend formula
\[ i_\F(B)(X) \;\;=\;\; \int^{n\in \F} X^n\times C(n] \]
where $B\in Set^\F$. The following calculation shows that the right hand square in the above diagram commutes:
\[ i_\F(Lan_\gamma(C)(X) \cong \int^{m\in \F} X^m\times Lan_\gamma(C)(m] \cong \]
\[ \cong \int^{m\in \F} X^m\times \sum_{n\in \o} (m]^n\otimes_n C(n] \cong\]
\[ \cong \sum_{n\in \o} X^n\otimes_n C(n] \]
The functor  $(-)^{sa}$, right adjoint to $i^{sa}$, is given by the formula
\[ F^{sa}(X) = \sum_{n\in\o} \left[\begin{array}{c} X \\ n  \end{array}\right]\otimes_n F(n] \]
for $F\in \End$. $(-)^{sa}$ associates to functors and natural transformations their 'semi-analytic parts'.

Note that both $\San$ and $\End$ are strict monoidal categories with tensor given by composition, and $i^{sa}$ is a strict monoidal functor. Thus its right adjoint $(-)^{sa}$ has a unique lax monoidal structure making the adjunction $i^{sa}\dashv (-)^{sa}$ a monoidal adjunction. This in turn gives us a monoidal monad $(\W,\eta,\mu)$ on $\San$.

We have (cf. \cite{Z}) a $2$-natural transformation $\cU$
\begin{center} \xext=1000 \yext=350
\begin{picture}(\xext,\yext)(\xoff,\yoff)
%1-cells
\putmorphism(0,250)(1,0)[\phantom{\MonCat}`\phantom{\bCat}`\Mon]{1000}{1}a
\putmorphism(0,150)(1,0)[\MonCat`\bCat`]{1000}{0}a
\putmorphism(0,50)(1,0)[\phantom{\MonCat}`\phantom{\bCat}`|-|]{1000}{1}b
%2-cell
\put(500,100){\makebox(300,100){$\Da\; \cU$}}
\end{picture}
\end{center}
where $\MonCat$ is the $2$-category of monoidal categories, lax monoidal functors, and monoidal transformations; $\Mon$ is the $2$-functor associating the monoids objects to monoidal categories, $|-|$ is the forgetful functor forgetting the monoidal structure, and $\cU$ is a $2$-natural transformation whose component at a monoidal category $M$ is the forgetful functor from monoids in $M$ to the underlying category of $M$:  $\cU_M: \Mon(M) \ra |M|$.

Applying $\cU$ to the monoidal adjunction and $i^{sa}\dashv (-)^{sa}$ and monoidal monad $\W$ we get an adjunction between categories of monoids and a monad on $\Mon(\San)$. The unnamed arrow in the above diagram is $\Mon(i^{sa})$. But the monoids in $\End$ and $\San$ are exactly monads and hence we get the left most adjunction $\P_f^m\dashv \Q_r^m$ together with the monad $(\Wb,\bar{\eta},\bar{\mu})$ on the category of semi-analytic monads.

On the other hand, on the categories $Set^\F$ and $Set^\S$ there are substitution tensors making $i_\F$ and $\i_\S$ monoidal equivalences and $\gamma^*$ and $Lan_\gamma$ monoidal adjunctions. Thus we can apply the 2-functor $\Mon$ to this adjunction and obtain an adjunction $\Mon( Lan_\gamma)\dashv \Mon(\gamma^*)$
as in the diagram

\begin{center} \xext=1900 \yext=900
\begin{picture}(\xext,\yext)(\xoff,\yoff)
%squares
\setsqparms[-1`0`0`-1;1600`750]
 \putsquare(150,50)[Set^\F`\Mon(Set^\F)\; =\; \FOp\hskip 25mm`Set^\S`\Mon(Set^\S)\; = \; \RegOp\hskip 25mm;\cU_{Set^\F}```\cU_{Set^\S}]

%vertical
 \putmorphism(1200,750)(0,-1)[\phantom{\Mon(Set^\F)}`\phantom{\Mon(Set^\S)}`]{700}{1}r
 \putmorphism(1100,750)(0,-1)[\phantom{\Mon(Set^\F)}`\phantom{\Mon(Set^\S)}`\Mon(Lan_\gamma)]{700}{-1}l
 \putmorphism(1800,750)(0,-1)[\phantom{\FOp}`\phantom{\RegOp}`\Q_r^o]{700}{1}r
 \putmorphism(1700,750)(0,-1)[\phantom{\FOp}`\phantom{\RegOp}`\P_f^o]{700}{-1}l
  \putmorphism(200,750)(0,-1)[\phantom{Set^\F}`\phantom{Set^\S}`\gamma^*]{700}{1}r
 \putmorphism(100,750)(0,-1)[\phantom{Set^\F}`\phantom{Set^\S}`Lan_\gamma]{700}{-1}l
 \end{picture}
\end{center}
The unnamed functor is $\Mon((-)^{sa})$. But monoids in $Set^\F$ and $Set^\S$ are (equivalent to) the categories of full and regular operads, respectively. The verification that the right most square commutes serially is left for the reader. We obtain

\begin{proposition}
The square (\ref{square_semi_an}) of categories and functors commutes up to an isomorphism.
\end{proposition}

{\em Proof.}  Both horizontal adjunctions in the square (\ref{square_semi_an}) are obtained from equivalent monoidal adjunctions. It remains to show that the identifications we obtained above are isomorphic to the functors $\cM^f_o$ and $\cM^r_o$, respectively. This is left for the reader. $\boxempty$

\vskip 2mm

There are free monads on finitary functors (cf. \cite{Barr}) and free semi-analytic monads on semi-analytic functors. The adjunctions $F\dashv U$ and $\widehat{F}\dashv \widehat{U}$ induce monads $\R$ and $\widehat{\R}$, respectively. $\widehat{\R}$ is the finitary version of what is called 'the monad for all monads' in \cite{Barr}. Putting this additional data to the above diagram and simplifying it at the same time we get a diagram

\begin{center} \xext=3000 \yext=1450
\begin{picture}(\xext,\yext)(\xoff,\yoff)
\setsqparms[0`0`0`0;1600`750]
 \putsquare(450,450)[\Mnd`\End`\SanMnd`\San;```]

%horizontal
\putmorphism(450,1250)(1,0)[\phantom{\Mnd}`\phantom{\End}`\widehat{F}]{1600}{-1}a
\putmorphism(450,1150)(1,0)[\phantom{\Mnd}`\phantom{\End}`\widehat{U}]{1600}{1}b

\putmorphism(450,500)(1,0)[\phantom{\AMnd}`\phantom{\San}`F]{1600}{-1}a
\putmorphism(450,400)(1,0)[\phantom{\AMnd}`\phantom{\San}`U]{1600}{1}b

%vertical
 \putmorphism(500,1150)(0,-1)[\phantom{\Mnd}`\phantom{2\bCat_{\times}}`\Q^m_r]{700}{1}r
 \putmorphism(400,1150)(0,-1)[\phantom{\Mnd}`\phantom{2\bCat_{\times}}`\P^m_f]{700}{-1}l
 \putmorphism(2100,1150)(0,-1)[\phantom{\End}`\phantom{\San}`(-)^{sa}]{700}{1}r
 \putmorphism(2000,1150)(0,-1)[\phantom{\End}`\phantom{\San}`i^{sa}]{700}{-1}l

%\Wb
  \put(450,50){\oval(100,100)[b]}
  \put(400,50){\line(0,1){300}}
  \put(500,50){\vector(0,1){300}}
  \put(0,150){${(\Wb,\bar{\eta},\bar{\mu})}$}

%\W
  \put(2050,50){\oval(100,100)[b]}
  \put(2000,50){\line(0,1){300}}
  \put(2100,50){\vector(0,1){300}}
  \put(2150,150){${(\W,\eta,\mu)}$}

%\M
  \put(2550,450){\oval(100,100)[r]}
  \put(2550,500){\line(-1,0){400}}
  \put(2550,400){\vector(-1,0){400}}
  \put(2400,550){${(\R,\eta,\mu)}$}

%\M_hat
  \put(2550,1200){\oval(100,100)[r]}
  \put(2550,1250){\line(-1,0){400}}
  \put(2550,1150){\vector(-1,0){400}}
  \put(2400,1300){${(\widehat{\R},\widehat{\eta},\widehat{\mu})}$}

 \end{picture}
\end{center}
In the above diagram the square of the left adjoints commutes. Thus, the square of the right adjoint commutes as well.
This shows in particular that the free monad on a semi-analytic functor is semi-analytic.

The monad $\Wb$ is a lift of a monad $\W$ to the category of $\R$-algebras $\SanMnd$ and, by \cite{Beck} we obtain

\begin{theorem}
The monad $\R$ for regular monads distributes over the monad  $\W$ for finitary functors, i.e. we have a distributive law
\[ \lambda: \R\W \lra \W\R \]
The category of algebras of the composed monad $\W\R$ on $\SanMnd$ is equivalent to the category $\Mnd$ of all finitary monads on $Set$. $\boxempty$
\end{theorem}

\vskip 2mm
{\em Remark.} We arrived at the above theorem with essentially no calculations.  We give below explicit formulas how to calculate
the values of some functors mentioned above and we shall also describe the coherence morphism on the monoidal monad $\W$. The coherence morphism $\varphi : \W \circ \W \lra \W(-\circ-)$  is a 'finite' description of the distributive law $\lambda$.

\vskip 2mm

First we describe the adjunction $i^{sa}\dashv (-)^{sa}$. We shall drop the inclusion $i^{sa}$ when possible.
Let $A\in \San$ and $G\in \End$ and $X$ be a set. The regular functor $A$ is given by the (functor of) its coefficients. Its value at $X$ is
\[ A(X)=\sum_{n\in \o} \left[\begin{array}{c} X \\ n  \end{array}\right]\otimes_n A_n\]
The value of $G^{sa}$ at $X$ is
\[ G^{sa}(X)= \sum_{n\in \o} \left[\begin{array}{c} X \\ n  \end{array}\right]\otimes_n G(n]\]
Thus
\[ \W(A)(X)=A^{sa}(X) =\sum_{n,m\in \o} \left[\begin{array}{c} X \\ n  \end{array}\right]\otimes_n(n]^m\otimes_m A_m  \]

The unit of the adjunction $i^{sa}\dashv (-)^{sa}$ at $X$
\[ (\eta_A)_X : A(X) \lra A^{sa}(X) \]
is given by
\[ [\vec{x},a] \mapsto [\vec{x},1_n,a]   \]
where $\vec{x}:(n]\ra X$ and $a\in A_n$.

The counit of the adjunction at $X$
\[ (\varepsilon_G)_X : \sum_{n\in\o} \left[\begin{array}{c} X \\ n  \end{array}\right] \otimes_n G(n]\lra G(X) \]
is given by
\[ [\vec{x},t] \mapsto G(\vec{x})(t)   \]
where $\vec{x}:(n]\ra X$ is an injection and $t\in G(n]$.

The multiplication in the monad $\W$
\[ (\mu_A)_X : \sum_{n,m,k\in \o} \left[\begin{array}{c} X \\ n  \end{array}\right]\otimes_n(n]^m\otimes_m(m]^k\otimes_k A_k  \lra \sum_{n,k\in \o} \left[\begin{array}{c} X \\ n  \end{array}\right]\otimes_n(n]^k\otimes_k A_k \]
is given by composition
\[ [\vec{x},g,f,a]\mapsto  [\vec{x},g\circ f,a]  \]
where $\vec{x}:(n]\ra X$, $g:(m]\ra (n]$, $f:(k]\ra (m]$, and $a\in A_k$.
This ends the definition of the monad $\W$.

Now we shall describe the monoidal structure on $\W$. If $B$ is another analytic functor, the $n$-coefficient of the composition $A\circ B$ is given by
\[ (A\circ B)_n =  \sum_{m,n_1,\ldots, n_m\in \o, \sum_{i=1}^m n_i =n} (S_n\times B_{n_1}\times\ldots \times B_{n_m}\times A_m)_{/\sim_n}  \]
where the equivalence relation $\sim_n$ is such that for $\sigma\in S_n$, $\sigma_i\in S_{n_i}$, $\tau\in S_m$, $b_i\in B_i$, for $i\in (m]$ and $a\in A_m$ we have
\[  \lk \sigma, \sigma_1\cdot b_1,\ldots, \sigma_m\cdot b_m, \tau\cdot a \rk \sim_n
\lk \sigma\circ (\sigma_{1},\ldots, \sigma_{m})\star\tau, b_{\tau(1)},\ldots, b_{\tau(m)}, a \rk  \]

The $n$-th coefficient of  $\W(A)\circ \W(A)$ is given by
\[ (\W(A)\circ \W(A))_n =  \sum_{m,n_i,k_i\in \o, \sum_{i=1}^m n_i =n} (S_n\times (n_1]^{k_1} \otimes_{k_1} A_{k_1}\times\ldots \times (n_m]^{k_m} \otimes_{k_m} A_{k_m}\times A_m)_{/\sim_n}  \]
and the  $n$-coefficient of  $\W(A^2)$ is given by
\[ (\W(A^2))_n =  \sum_{m,k,k_i\in \o, \sum_{i=1}^m k_i =k} (n]^k\times A_{k_1}\times\ldots \times A_{k_m} \times A_m \]
The coherence morphism $\varphi$ for $\W$ at the $n$-th coefficient of the functor $A$ is
\[ \varphi_n: (\W(A)\circ \W(A))_n  \lra (\W(A^2))_n  \]
is given by
\[  \lk \sigma, [\sigma_1,a_1]\ldots [\sigma_m,a_m], \tau, a \rk \mapsto
\lk \sigma\circ(\sigma_{1}\ldots \sigma_{m})\star\tau, a_{\tau(1)}\ldots a_{\tau(m)}, a\rk   \]
Note that this map is well defined at the level of equivalence classes.

%\newpage
\section{Equational theories vs regular operads}\label{sec_eqth}

In this section we study the relations between regular equational theories  and regular operads. We shall show that the square
\begin{equation}\label{square-eo}\end{equation}
\begin{center} \xext=800 \yext=400
\begin{picture}(\xext,\yext)(\xoff,\yoff)
\setsqparms[1`-1`-1`1;800`500]
\putsquare(0,50)[\FOp`\ET`\RegOp`\RegET;\cE_o^f`\P_f^o`\P_f^e`\cE_o^r]
 \end{picture}
\end{center}
commutes up to an isomorphism, with $\P_f^e$ being an inclusion and both horizontal functors being equivalences of categories. $\P_f^o`$ was defined in Section \ref{pres_of_cats_of_algs} and $\cE_o^f$ was defined in Section \ref{full_operads}. We shall define $\cE_o^r$.

\subsection*{The functor $\cE_o^r : \RegOp \ra \ET$}

Let $\cO$ be a regular  operad. We define an equational theory $\cE_o^r(\cO)=(L,A)$. As the set of $n$-ary function symbols we put $L_n=\cO_n$ for $n\in \o$. The set of axioms $A$ contains the following three kinds of equations in context:
\begin{enumerate}
  \item {\em unit}: $\iota(x_1)=x_1 :  \vec{x}^1$ where $\iota\in \cO_1$ is the unit of the operad $\cO$;
  \item {\em action}: $ a(x_{\tau(1)}, \ldots , x_{\tau(m)})=(\tau\cdot a)(x_1, \ldots , x_n) : \vec{x}^n$ for all $a\in \cO_m$ and surjections $\tau:(m] \ra (n]$;
   \item {\em multiplication}: $a(a_1(x_1, \ldots , x_{k_1}), \ldots ,a_m(x_{k-k_m+1}, \ldots , x_{k})) = ((a_1,\ldots, a_m)\ast a)(x_1, \ldots x_k) : \vec{x}^k$ where $a\in \cO_m$, $a_i\in \cO_{k_i}$ for $i\in 1,\ldots, m$, $k=\sum_{i=1}^mk_i$;
\end{enumerate}
Clearly, all equations are regular and hence the theory $\cE_o^r(\cO)$ is regular.

Suppose that $h:\cO\ra \cO'$ is a morphism of regular operads. We define the interpretation
$$\cE_o^r(h): \cE_o^r(\cO)\lra \cE_o^r(\cO')$$ For $a\in \cO_n$
we put
\[ \cE_o^r(h)(a) = h(a)(x_1,\ldots, x_n):\vec{x}^n\]
for $n\in \o$. Clearly, $\cE_o^r(h)$ is a regular interpretation.

\begin{proposition}\label{commuation_et_lt}
The square (\ref{square-eo}) commutes up to a natural isomorphism.
\end{proposition}
{\em Proof.}  Let $\cO$ be a regular operad. We define an interpretation of equational theories
\[ I_\cO : \cE^r_o(\cO) \lra \cE^f_o \P^o_f(\cO) \]
by
\[ \cO_n\ni a  \mapsto [[1_n,a](x_1,\ldots x_n):\vec{x}^n] \]
where $a\in \cO_n$ and $[1_n,a]$ is an $n$-ary operation symbol of the theory $\cE^f_o \P^o_f(\cO)$.

We need to verify that $I$ is a well defined natural isomorphism. First we need to verify that $I_\cO$ preserves axioms. The unit axiom is obvious.
To prove the action axioms, we fix $a\in \cO_m$ and $\tau: (m]\ra (n]$ and we calculate using the theory $\cE^f_o \P^o_f(\cO)$
\[ I_\cO(a(x_{\tau(1)},\ldots ,x_{\tau(m)})= [1_m,a](x_{\tau(1)},\ldots ,x_{\tau(m)})= \]
\[ = \tau\cdot[1_m,a](x_1,\ldots ,x_n)= [1_m,\tau\cdot a](x_1,\ldots ,x_n)=  \]
\[ = I_\cO(\tau\cdot a(x_1,\ldots ,x_n))\]
The calculations for the composition axioms are similar. The naturality of $I_\cO$ is left for the reader.

We shall show that $I_\cO$ is an isomorphism of theories. To this end we define an inverse interpretation $$J_\cO:  \cE^f_o \P^o_f(\cO)\lra \cE^r_o(\cO)$$
given by
\[ \I(k,n)\otimes_k \cO_k \ni [f,a] \mapsto a(x_{f(1)},\ldots, x_{f(k)}):\vec{x}^n \]
Again we need to verify that $J_\cO$ preserves the axioms and again we shall verify the action axiom only. Fix $\phi:[n)\ra [m)$ and $[f,a]\in \I(k,n)\otimes \cO_k$. Let $f'\circ\phi'$ be an epi-mono factorization in $\F$ with $f':[n')\ra [m)$. Using the theory $\cE^r_o(\cO)$, we have
\[ J_\cO((\phi\cdot[f,a])(x_1,\ldots ,x_n)) = [f',\phi'\cdot a](x_1,\ldots ,x_k) = \]
\[ = (\phi'\cdot a)(x_{f'(1),\ldots,x_{f'(n')}}) = a(x_{f'\circ\phi'(1),\ldots,x_{f'\circ\phi'(n')}}) =\]
\[ = a(x_{\phi\circ f(1),\ldots,x_{\phi\circ f(n')}}) = J_\cO([f, a](x_{\phi(1),\ldots,x_{\phi(n')}}))  \]

Finally, we need to verify that $I_\cO$ and $J_\cO$ are mutually inverse one to the other. The composition $J_\cO I_\cO$ sends operation $a\in \cO_n$ to the term $a(x_1,\ldots ,x_n):\vec{x}^n$ so it is the identity.
For an operation $[f,a]\in \I(k,n)\otimes_k \cO_k$ in theory $\cE^f_o \P^o_f(\cO)$ we have
\[ I_\cO J_\cO([f,a])= [1_n,a](x_{f(1)},\ldots,x_{f(n)}):\vec{x}^n = \]
\[ = f\cdot [1_n,a](x_1,\ldots,x_k):\vec{x}^k = [f,a](x_1,\ldots,x_k):\vec{x}^k\]
i.e. the composition $I_\cO J_\cO$ is the identity as well.
 $\boxempty$

Finally, we have

\begin{theorem}
The functor $\cE^r_o : \RegOp\lra \RegET$ is an equivalence of categories.
\end{theorem}

{\em Proof.} From Proposition \ref{commuation_et_lt} it follows that $\cE^r_o$ is faithful. To see that it is full consider a regular interpretation
$I:\cE^r_o(\cO)\ra \cE^r_o(\cO')$. We shall define a morphism of regular operads $h:\cO\ra\cO'$ such that $\cE^r_o(h)=I$. For $f\in \cO_n$, $I(f)(\vec{x}^n):\vec{x}^n$ is a regular term. But in theories in the image of $\cE^r_o$ every regular term is equal to a simple term, i.e. there is an $f'\in \cO'_n$ such that $f'(\vec{x}^n)=I(f)(\vec{x}^n):\vec{x}^n$ is a theorem in the theory $\cE^r_o(\cO')$. Thus we put $h(f)=f'$. It is left for the reader to verify that $h$ has the required property.

To see that $\cE^r_o$ is essentially surjective let us fix a regular theory $T=(L,A)$. Then the regular terms in $T$ form a regular operad called $T^{ro}$. The unit is the term $x_1:x_1$.  The composition is defined by the substitution (making sure that we make the variable disjoint in different substituted term, via $\alpha$-conversion). The action of a surjection $\phi:(n]\ra(m]$ on a regular term $t(x_1,\ldots,x_n):\vec{x}^n$  is again a regular term $t(x_{\phi(1)},\ldots,x_{\phi(n)}:\vec{x}^n$. Again it is a matter of a routine verification that  $\cE^r_o(T^{ro})\cong T$.
$\boxempty$

\vskip 2mm
\subsection*{Examples}
\begin{enumerate}
    \item Let ${\mathds{1}}$ be the terminal equational theory. It has one constant, say $e$, and can be axiomatized by a single axiom: $v_1=e:\vec{v}^1$. As a Lawvere theory it is the category that has exactly one morphism between any two objects. It is best seen at the level of Lawvere theories. Both theories  $\mathds{1}^a$ and the theory of commutative monoids are linear-regular and have exactly one analytic morphism $a:n\ra 1$, for any $n$.

  \item   The functor $\Q^e_r: \ET\ra \RegET$ is a right adjoint, it preserves the terminal object. Hence $\mathds{1}^r$, the regular part  ${\mathds{1}}$, is the terminal regular theory. It is the theory of suplattices.

  The embedding of the regular theories into all equational theories has a right adjoint $(-)^{sr}$, as well. The value of the functor $(-)^{sr}$ on the terminal equational theory  ${\mathds{1}}$ is the terminal regular theory, i.e. the theory of suplattices:
      \[ (x\vee y)\vee z = x\vee (y\vee z), \;\;\; x\vee y = y\vee x, \;\;\; x\vee \perp =x = \perp \vee x,\;\;\; (x\vee x) = x \]

  \item  The theory ${\mathds{1}}$ has a proper subcategory,  in which $0\not\cong 1$. It has no function symbols, and can be axiomatized by a single axiom: $v_1=v^2:\vec{v}^2$.  The regular part of this theory is the theory of suplattices without a bottom element.
\end{enumerate}

\section{Cartesian and weakly cartesian monads}\label{sec_cart}

In this section we shall investigate two (strict monoidal) subcategories of $\San$ and their categories of monoids. The category of (weakly) cartesian functors and (weakly) cartesian natural transformations will be denoted by $\Cart$ ($\wCart$). The corresponding categories of monoids: the category of (weakly) cartesian monads will be denoted by $\CartMnd$ ($\wCartMnd$). Thus we have embeddings full on isomorphisms
\[ \Cart\lra \wCart \lra \San \]
which are strict monoidal and induce embeddings of categories of monoids
\[ \CartMnd\lra \wCartMnd \lra \SanMnd \]
The characterizations of the subcategories of equational theories $\ET$ and of Lawvere theories $\LT$  corresponding to $\CartMnd$ and $\wCartMnd$ are a bit technical and we are not going to describe it in detail here. Clearly, the objects are some regular theories satisfying additional conditions and similarly for morphisms.
We shall content ourselves with a description of subcategories of $Set^\S$ whose essential images are  $\wCart$ and $\Cart$, respectively. Note, however, that if $(T,\eta,\mu)$ is a semi-analytic monad such that the functor part $T$ is the left Kan extension of a functor $R:\S\ra Set$ then $R$ is the functor of all regular operations in the equational theory corresponding to the monad $T$. Thus our description will in fact provide a description of the equational theories corresponding to monads in $\wCart$ and $\Cart$.

{\em Remarks.} If in a weak pullback square
\begin{center} \xext=500 \yext=500
\begin{picture}(\xext,\yext)(\xoff,\yoff)
\setsqparms[1`1`1`1;500`400]
 \putsquare(0,50)[P`A`B`C;m`h`f`g]
 \end{picture}
\end{center}
the morphism $m$ is mono then the square is a pullback. Since $m$ is mono iff the square
\begin{center} \xext=500 \yext=500
\begin{picture}(\xext,\yext)(\xoff,\yoff)
\setsqparms[1`1`1`1;500`400]
 \putsquare(0,50)[P`P`P`A;1_P`1_P`m`m]
 \end{picture}
\end{center}
is a pullback, it follows that if a functor weakly preserves pullbacks it does preserve monos as well.
Recall the description of functor $\hat{(-)}: Set^\S \lra \End$ from Section \ref{monads}, the Kan extension along $i_\S:\S\ra Set$. We begin with the following observation

\begin{lemma} \label{equal} Let $A:\S\ra Set$ be a functor, $f:X\ra Y$ a function, $[\vec{x}:(n]\ra X,a\in A_n], [\vec{x}':(n']\ra X,a'\in A_{n'}]\in \hat{A}(X)$. If
\[ \hat{A}([\vec{x},a]) = \hat{A}([\vec{x}',a']) \]
then there are $m\in \o$, surjections  $g:n\ra m$, $g':n'\ra m$, and an injection $\vec{y}:m\ra Y$ as in the diagram
\begin{center} \xext=800 \yext=700
\begin{picture}(\xext,\yext)(\xoff,\yoff)
\setsqparms[2`3`1`2;800`600]
 \putsquare(0,50)[(n]`X`(m]`Y;\vec{x}`g`f`\vec{y}]
  \put(200,400){$(n']$}
  \put(350,450){\vector(3,1){400}}
  \put(350,450){\vector(3,1){40}}
  \put(500,380){$\vec{x}'$}
  \put(200,350){\vector(-1,-2){120}}
  \put(200,350){\vector(-1,-2){110}}
  \put(180,200){$g'$}
 \end{picture}
\end{center}
such that $\vec{y}\circ g = f\circ \vec{x}$ $\vec{y}\circ g' = f\circ \vec{x}'$, and $A(g)(a)=A(g')(a')$.
In particular
\[ \hat{A}([\vec{x},a]) =[\vec{y},A(g)(a)]=[\vec{y},A(g')(a')]= \hat{A}([\vec{x}',a']) \]
\end{lemma}

{\em Proof.} Exercise. $\boxempty$

The following two Propositions identify the subcategory of $Set^\S$ whose essential image in $\End$ is $\wCart$.

\begin{proposition}\label{im-wCart} Let $A:\S\ra Set$ be a functor. The functor $\hat{A} : Set \ra Set$ weakly preserves pullbacks iff the functor $A$ satisfies the following condition (WPB):

for any pullback of surjections in $\F$
\begin{center} \xext=500 \yext=500
\begin{picture}(\xext,\yext)(\xoff,\yoff)
\setsqparms[1`1`1`1;500`400]
 \putsquare(0,50)[(p]`(m]`(n]`{(k]};g_2`f_2`f_1`g_1]
 \end{picture}
\end{center}
and a pair of elements $a\in A(m]$, $b\in A(n]$ such that $A(f_1)(a)=A(g_1)(b)$ there is an injection $h: (q]\ra (p]$ and an element $c\in A(q]$ such that
$f_2\circ h$ and $g_2\circ h$ are surjections and
\[ A(g_2\circ h)(c) =a, \;\;\; {\rm and}\;\;\;  A(f_2\circ h)(c) =b\]
\end{proposition}

{\em Proof.} $\Ra$ First we verify that the condition (WPB) is necessary. Let us fix the square and elements $a$ and $b$ as in the above condition. Then we have elements
\[ [1_{(m]},a]\in \left[\begin{array}{c} (m] \\ m  \end{array}\right]\otimes_m A(m)\subseteq \sum_{k\in\o} \left[\begin{array}{c} (m] \\ k  \end{array}\right]\otimes_k A(k) =\hat{A}(X)  \]
and
\[ [1_{(n]},b]\in \left[\begin{array}{c} (n] \\ n  \end{array}\right]\otimes_n A(n)\subseteq \hat{A}(X)  \]
such that
\[ \hat{A}(f_1)([1_{(m]},a]) =\hat{A}(g_1)([1_{(n]},b]) \]
As $f_1$, $g_1$ are surjections this means that
\[ A(f_1)(a)=A(g_1)(b) \]
As $\hat{A}$ weakly preserves pullbacks there is an $[h,c]\in \hat{A}(p]$ such that
\[ \hat{A}(f_2)([h,c]) = [1_{(n]},b], \;\;\;{\rm and}\;\;\; \hat{A}(g_2)([h,c]) = [1_{(m]},a] \]
Thus for some $h:(q]\ra (p]$, $c\in A(q)$ we have $[h,c]\in  \left[\begin{array}{c} (p] \\ q  \end{array}\right]\otimes_q A(q)$
and using the diagram
\begin{center} \xext=800 \yext=900
\begin{picture}(\xext,\yext)(\xoff,\yoff)
\setsqparms[1`1`1`1;400`400]
 \putsquare(400,50)[(p]`(m]`(n]`{(k]};g_2`f_2`f_1`g_1]
 \setsqparms[1`1`0`0;800`800]
 \putsquare(0,50)[(q]`(m']`(n']`;g_3`f_3``]
  \put(100,780){\vector(1,-1){260}}
  \put(260,650){$h$}
  \putmorphism(0,50)(1,0)[\phantom{(n']}`\phantom{(n]}`\vec{y}]{400}{1}b
  \putmorphism(800,850)(0,-1)[\phantom{(n']}`\phantom{(n]}`\vec{x}]{400}{1}r
 \end{picture}
\end{center}
we have
\[ [\vec{y},A(f_3)(c)]= [1_{(n]},b]\;\;\;{\rm and}\;\;\; [\vec{x},A(g_3)(c)]= [1_{(m]},a]  \]
This means that $n=n'$, $\vec{y}$ is a bijection, $\vec{y}\circ f_3= f_2\circ h$ is a surjection and $b=A(\vec{y}\circ f_3)=A(f_2\circ h)$. Similarly, we get that $m=m'$, $g_2\circ h$ is a surjection and $a=A(g_2\circ h)$. Thus the condition (WPB) is necessary.

$\La$ To show that the condition (WPB) is sufficient we suppose that a functor $A:\S\ra Set$ satisfies (WPB) and we shall show that $\hat{A}: Set \ra Set$ weakly preserves pullbacks.

Let us fix a pullback in $Set$ and
\begin{equation}\label{pb}\end{equation}
\begin{center} \xext=500 \yext=350
\begin{picture}(\xext,\yext)(\xoff,\yoff)
\setsqparms[1`1`1`1;500`400]
 \putsquare(0,50)[R`Y`Z`X;g_2`f_2`f_1`g_1]
 \end{picture}
\end{center}
and elements $[\vec{z},b]\in \hat{A}(Z)$ and $[\vec{y},a]\in \hat{A}(Y)$ such that
\[  \hat{A}(f_1)([\vec{y},a])=\hat{A}(g_1)([\vec{z},b]) \]
i.e. in the diagram
\begin{equation}\label{diag7cos}\end{equation}
\begin{center} \xext=1650 \yext=1200
\begin{picture}(\xext,\yext)(\xoff,\yoff)
\putmorphism(0,1250)(1,0)[(r]`\phantom{(q]}`h]{400}{1}a
\setsqparms[1`1`0`0;1200`1200]
\putsquare(400,50)[(q]`(m]`(n]`;g_4`f_4``]
\setsqparms[1`1`1`1;400`400]
\putsquare(800,450)[R`Y`Z`X;g_2`f_2`f_1`g_1]

\putmorphism(400,50)(1,0)[\phantom{(n]}`(k]`g_3]{800}{1}b
\putmorphism(1600,1250)(0,-1)[\phantom{(m]}`(k]`g_3]{800}{1}r

\put(500,1180){\vector(1,-1){260}}
\put(650,1050){$v$}

\put(1500,1180){\vector(-1,-1){260}}
\put(1300,1050){$\vec{y}$}

\put(500,110){\vector(1,1){260}}
\put(600,260){$\vec{z}$}

\put(1200,150){\vector(0,1){240}}
\put(1110,230){$\vec{z}'$}

\put(1500,450){\vector(-1,0){260}}
\put(1400,500){$\vec{y}'$}

\put(1300,80){\vector(1,1){260}}
\put(1410,130){$\sigma$}
 \end{picture}
\end{center}
we have
\[  [\vec{y'},A(f_3)(a)]=[\vec{z'},A(g_3)(b)] \]
i.e. there is a permutation $\sigma\in S_k$ such that
\[ \vec{y}\circ \sigma = \vec{z'}\;\;\;{\rm and}\;\;\; A(f_3)(a)=A(\sigma\circ g_3)(b) \]
Still in the above diagram we can form a pullback in $\F$
\begin{center} \xext=500 \yext=500
\begin{picture}(\xext,\yext)(\xoff,\yoff)
\setsqparms[1`1`1`1;500`400]
 \putsquare(0,50)[(q]`(m]`(n]`{(k]};g_4`f_4`f_3`\sigma\circ g_3]
 \end{picture}
\end{center}
and we have a morphism $v: (q]\ra R$ into the pullback (\ref{pb}) such that
\[ f_2\circ v = \vec{z}\circ f_4,\;\;\; g_2\circ v = \vec{y}\circ g_4 \]
Thus by assumption there is an injection $h:(r]\ra (q]$ and a $c\in A(q)$ such that both $f_4\circ h$ and $g_4\circ h$ are surjections
and
\[ A(g_4\circ h)(c)=a,\;\; \; A(f_4\circ h)(c)=b\]
From diagram (\ref{diag7cos}) we see that $f_4\circ h$ is a surjection, $\vec{z}$ is an injection and
$f_2\circ v \circ h =\vec{z}\circ f_4 \circ h$. Thus
\[ \hat{A}(f_2)(\hat{A}(v)([h,c])= [\vec{z}, A(f_4\circ h)(c)=[\vec{z},b] \]
Similarly
\[ \hat{A}(g_2)(\hat{A}(v)([h,c])=[\vec{y},a] \]
Thus $\hat{A}(v)([h,c])\in \hat{A}(P)$ is the sought element in the weak pullback
\begin{center} \xext=1600 \yext=600
\begin{picture}(\xext,\yext)(\xoff,\yoff)
\setsqparms[1`1`1`1;600`400]
 \putsquare(1000,50)[\hat{A}(P)`\hat{A}(Y)`\hat{A}(Z)`\hat{A}(X);\hat{A}(g_2)`\hat{A}(f_2)`\hat{A}(f_1)`\hat{A}(g_1)]
  \putmorphism(0,450)(1,0)[{[h,c]\in\hat{A}(q]}`\phantom{\hat{A}(P)}`\hat{A}(v)]{1000}{1}a
 \end{picture}
\end{center}
$\boxempty$

\begin{proposition} \label{wpb}Let $\tau : A\ra B$ be a natural transformation. Then $\tau$ is weakly cartesian iff  $\hat{\tau}: \hat{A} \ra \hat{B}$
is weakly cartesian.
\end{proposition}

{\em Proof.} Assume that $\tau:A\ra B$ is a weakly cartesian natural transformation, $f:X\ra Y$ a function. We shall show that the square
\begin{center} \xext=600 \yext=550
\begin{picture}(\xext,\yext)(\xoff,\yoff)
\setsqparms[1`1`1`1;600`400]
 \putsquare(0,50)[\hat{A}(X)`\hat{B}(X)`\hat{A}(Y)`\hat{B}(Y);\hat{\tau}_X`\hat{A}(f)`\hat{B}(f)`\hat{\tau}_Y]
 \end{picture}
\end{center}
is a weak pullback. Let us fix elements
\[ [\vec{y},a] \in \left[\begin{array}{c} Y \\ n  \end{array}\right]\otimes_n A(n) \subseteq \hat{A}(Y),\;\;\;
{\rm and}\;\;\; [\vec{x},b] \in \left[\begin{array}{c} X \\ m  \end{array}\right]\otimes_m B(m) \subseteq \hat{B}(X) \]
such that
\[ \hat{B}(f)([\vec{x},b]) = \hat{\tau}([\vec{y},a])\;\; (= [\vec{y},\tau_n(a)]) \]
By Lemma \ref{equal} we have a surjection $f'$
\begin{center} \xext=1600 \yext=550
\begin{picture}(\xext,\yext)(\xoff,\yoff)
\setsqparms[1`1`1`1;600`400]
 \putsquare(0,50)[(m]`X`(n]`Y;\vec{x}`f'`f`\vec{y}]
 \end{picture}
\end{center}
such that $B(f')(b)=\tau_n(a)$. Since $\tau$ is weakly cartesian, there is a $c\in A(m)$ such that
\[ \tau_m(c)=b, \;\;\;{\rm and}\;\;\; A(f')(c)=a \]
Then \[ [\vec{x},c]\in \left[\begin{array}{c} X \\ m  \end{array}\right]\otimes_m A(m) \]
and moreover
\[ \hat{A}(f)([\vec{x},c])= [\vec{y},A(f')(c)]=[\vec{y},a]\]
and
\[ \hat{\tau}_X([\vec{x},c]) = [\vec{x},\tau_m(c)]= [\vec{x},b] \]

To prove the converse let us assume that $\hat{\tau}:\hat{A}\ra \hat{B}$ is a weakly cartesian natural transformation and that $f:(m]\ra (n]$ is a surjection in $\S$. We need to show that the square
\begin{center} \xext=1600 \yext=550
\begin{picture}(\xext,\yext)(\xoff,\yoff)
\setsqparms[1`1`1`1;600`400]
 \putsquare(0,50)[A(m]`B(m]`A(n]`{B(n]};\tau_m`A(f)`B(f)`\tau_n]
 \end{picture}
\end{center}
is weakly cartesian. Fix $a\in A(n)$ and $b\in B(m)$ such that $\tau_n(a)=B(f)(b)$. Consider a weak pullback in $Set$
\begin{center} \xext=1600 \yext=550
\begin{picture}(\xext,\yext)(\xoff,\yoff)
\setsqparms[1`1`1`1;600`400]
 \putsquare(0,50)[\hat{A}(m]`\hat{B}(m]`\hat{A}(n]`{\hat{B}(n]};\hat{\tau}_m`\hat{A}(f)`\hat{B}(f)`\hat{\tau}_n]
 \end{picture}
\end{center}
We have elements
\[ [1_{(n]},a] \in \left[\begin{array}{c} (n] \\ n  \end{array}\right]\otimes_n A(n) \subseteq \hat{A}(n],\;\;\;
{\rm and}\;\;\; [1_{(m]},b] \in \left[\begin{array}{c} (m] \\ m  \end{array}\right]\otimes_m B(m) \subseteq \hat{B}(m] \]
such that
\[ \hat{\tau}_n([1_{(n]},a])= [1_{(n]},\tau_n(a)] = [1_{(n]},B(f)(b)]=\hat{B}(f)([1_{(m]},b]) \]
Thus, by our assumption, for some $k\in\o$ there is an $[\vec{x},c]\in \left[\begin{array}{c} (m] \\ k  \end{array}\right]\otimes_k A(k)$ such that
\[ \hat{A}(f)([\vec{x},c])=[1_{(n]},a],\;\;\;{\rm and}\;\;\; \hat{\tau}_m([\vec{x},c])=[1_{(m]},b] \]
Thus $(\vec{x},\tau_m(c))\sim (1_{(m]},b)$ and we have that $k=m$, $\vec{x}$ is a bijection, $(\vec{x},c)\sim (1_{(m]},A(\vec{x}(c))$. Hence we also have
$\tau_m(A(\vec{x})(c))=b$. Moreover as the square
\begin{center} \xext=500 \yext=500
\begin{picture}(\xext,\yext)(\xoff,\yoff)
\setsqparms[1`1`1`1;500`400]
 \putsquare(0,50)[(m]`(m]`(n]`{(n]};\vec{x}`f\circ\vec{x}`f`1_{(n]}]
 \end{picture}
\end{center}
commutes, we have
\[ A(f)(A(\vec{x})(c)) =A(f\circ\vec{x})(c)=a \]
Thus $A(\vec{x})(c)$ is the element sought for $a$ and $b$.  Since $f$, $a$ and $b$ were arbitrary, $\tau$ is weakly cartesian. $\boxempty$

The final two Propositions identify the subcategory of $Set^\S$ whose essential image  in $\End$ is $\Cart$.

\begin{proposition} Let $A:\S\ra Set$ be a functor. The functor $\hat{A} : Set \ra Set$ preserves pullbacks iff the functor $A$ satisfies the condition (WPB) from the Proposition \ref{im-wCart}, and additionally satisfies the following condition (PB):

suppose that the square
\begin{center} \xext=500 \yext=500
\begin{picture}(\xext,\yext)(\xoff,\yoff)
\setsqparms[1`1`1`1;500`400]
 \putsquare(0,50)[(p]`(m]`(n]`{(k]};g_2`f_2`f_1`g_1]
 \end{picture}
\end{center}
is a pullback of surjections in $\F$, $\vec{x}: (q]\ra (p]$, $\vec{x}':(q']\ra (p]$ are two injections, $c\in A(q]$, $c'\in A(q']$  are two elements so that
the functions
\[ f_2\circ \vec{x},\;\;\; f_2\circ \vec{x}',\;\;\; g_2\circ \vec{x},\;\;\; g_2\circ \vec{x}',\]
are surjections and
\[ A(f_2\circ \vec{x})(c) =   A(f_2\circ \vec{x}')(c'),\;\;\;\; A(g_2\circ \vec{x})(c) =   A(g_2\circ \vec{x}')(c').\]
Then $q=q'$ and there is $\sigma \in S_q$ such that
\[ \vec{x}'\circ \sigma = \vec{c},\;\;\;\; A(\sigma)(c)=c'. \]
\end{proposition}
{\em Proof.} Assume that $A$ satisfies (WPB) and (PB). Thus, by Proposition \ref{wpb}, $\hat{A}$ weakly preserves pullbacks.
Let
\begin{center} \xext=500 \yext=500
\begin{picture}(\xext,\yext)(\xoff,\yoff)
\setsqparms[1`1`1`1;500`400]
 \putsquare(0,50)[P`Y`Z`X;g_2`f_2`f_1`g_1]
 \end{picture}
\end{center}
be a pullback in $Set$. We shall show that the square
\begin{center} \xext=500 \yext=700
\begin{picture}(\xext,\yext)(\xoff,\yoff)
\setsqparms[1`1`1`1;700`500]
 \putsquare(0,50)[\hat{A}(P)`\hat{A}(Y)`\hat{A}(Z)`\hat{A}(X);\hat{A}(g_2)`\hat{A}(f_2)`\hat{A}(f_1)`\hat{A}(g_1)]
 \end{picture}
\end{center}
is a pullback, i.e. it satisfies the uniqueness condition. So let
\[ [h,a]\in \left[\begin{array}{c} P \\ q  \end{array}\right]\otimes_q A(q),\;\;\;{\rm and}\;\;\;  [h',a']\in \left[\begin{array}{c} P \\ q'  \end{array}\right]\otimes_{q'} A(q')  \]
be such that
\[ \hat{A}(f_2)([h,a])=\hat{A}(f_2)([h',a']),\;\;\;{\rm and}\;\;\; \hat{A}(g_2)([h,a])=\hat{A}(g_2)([h',a']) \]
This implies the equalities of images of functions
\[im(f_2\circ h)=im(f_2\circ h'),\;\;  im(g_2\circ h)=im(g_2\circ h'),\;\; im(g_1\circ f_2\circ h)=im(f_2\circ h')\]
All these sets are finite, say, having $n$, $m$, and  $k$ elements, respectively. Thus we can form a commuting diagram
\begin{center} \xext=1650 \yext=1500
\begin{picture}(\xext,\yext)(\xoff,\yoff)
\setsqparms[1`1`1`1;700`600]
\putsquare(950,50)[P`Y`Z`X;g_2`f_2`f_1`g_1]

\setsqparms[1`1`0`0;700`600]
\putsquare(400,450)[(r]`(m]`(n]`{(k]};g_4`f_4``]

%back to front
\put(500,980){\vector(3,-2){380}}
\put(1200,980){\vector(3,-2){380}}
\put(500,380){\vector(3,-2){380}}
\put(1200,380){\vector(3,-2){380}}

%back
\put(500,450){\line(1,0){420}}
\put(980,450){\vector(1,0){50}}

\put(1100,680){\line(0,1){300}}
\put(1100,620){\vector(0,-1){100}}

\put(980,800){$f_3$}
\put(600,500){$g_3$}

\putmorphism(0,1050)(1,0)[(q]`\phantom{(r]}`h]{400}{1}b
\putmorphism(400,1450)(0,-1)[(q']`\phantom{(r]}`h']{400}{1}l

\put(580,800){$v$}
\end{picture}
\end{center}
where the back square is a pullback of surjections (in $\F$) and all arrows from the back to the front are injections (for simplicity: inclusions). By definition the functions
\[ f_4\circ h,\;\; f_4\circ h',\;\; g_4\circ h,\;\; f_4\circ h' \]
are surjections. Since $\hat{A}$ preserves monos we also have
\[ A(f_4\circ h)(a)= A(f_4\circ h')(a'), \;\;\; A(g_4\circ h)(a)= A(g_4\circ h')(a') \]
Thus from the condition (PB) it follows that $q=q'$ and there is a $\sigma \in S_q$ such that $h'\circ \sigma = h$ and $A(\sigma)(a)=a'$, i.e.
$[h,a]=[h',a']$, as required.

To prove the converse we assume now that $\hat{A}$ preserves pullbacks and we fix a pullback
\begin{center} \xext=500 \yext=500
\begin{picture}(\xext,\yext)(\xoff,\yoff)
\setsqparms[1`1`1`1;500`400]
 \putsquare(0,50)[(r]`(m]`(n]`{(k]};g_2`f_2`f_1`g_1]
 \end{picture}
\end{center}
in $\F$ of surjections. As $\hat{A}$ weakly preserves pullbacks, the above pullback is sent by $\hat{A}$ to a weak pullback in $Set$. Then, using Lemma \ref{equal} it is easy to see that the condition (PB) expresses the fact that $\hat{A}$ sends the above square to a pullback in $Set$. $\boxempty$

\begin{proposition} Let $\tau : A\ra B$ be a natural transformation. Then $\tau$ is cartesian iff  $\hat{\tau}: \hat{A} \ra \hat{B}$
is cartesian.
\end{proposition}

{\em Proof.} First assume that $\tau : A\ra B$ is a cartesian natural transformation. Fix a function $f:X\ra Y$. By Proposition \ref{wpb} the square
\begin{center} \xext=600 \yext=550
\begin{picture}(\xext,\yext)(\xoff,\yoff)
\setsqparms[1`1`1`1;600`400]
 \putsquare(0,50)[\hat{A}(X)`\hat{B}(X)`\hat{A}(Y)`\hat{B}(Y);\hat{\tau}_X`\hat{A}(f)`\hat{B}(f)`\hat{\tau}_Y]
 \end{picture}
\end{center}
is a weak pullback. We shall show that it also satisfies uniqueness property. Let us fix elements
\[ [\vec{x},a] \in \left[\begin{array}{c} X \\ n  \end{array}\right]\otimes_n A(n) \subseteq \hat{A}(X),\;\;\;
{\rm and}\;\;\; [\vec{x}',a'] \in \left[\begin{array}{c} X \\ n'  \end{array}\right]\otimes_{n'} A(n') \subseteq \hat{A}(X) \]
such that
\begin{equation}\label{equ_pb}
\hat{\tau}([\vec{x},a])= \hat{\tau}([\vec{x}',a']),\;\;\;{\rm and}\;\;\; \hat{A}(f)([\vec{x},a]) =\hat{A}(f)([\vec{x}',a'])
\end{equation}
The first equality means that $[\vec{x},\tau_n(a)]= [\vec{x}',\tau_n(a')]$. Hence $n=n'$ and there is a $\sigma\in S_n$ such that
\[ \vec{x}'\circ \sigma = \vec{x},\;\;\; {\rm and}\;\;\; \tau_n(a')=B(\sigma)(\tau_n(a))= \tau_n(A(\sigma)(a)) \]
By Lemma \ref{equal} and the second equality in (\ref{equ_pb}) there are surjections $g$ and $g'$ and an injection $\vec{y}$ as in the diagram
\begin{center} \xext=800 \yext=700
\begin{picture}(\xext,\yext)(\xoff,\yoff)
\setsqparms[2`3`1`2;800`600]
 \putsquare(0,50)[(n]`X`(m]`Y;\vec{x}`g`f`\vec{y}]
  \put(200,400){$(n']$}
  \put(350,450){\vector(3,1){400}}
  \put(350,450){\vector(3,1){40}}
  \put(500,380){$\vec{x}'$}
  \put(200,350){\vector(-1,-2){120}}
  \put(200,350){\vector(-1,-2){110}}
  \put(180,200){$g'$}
  \put(80,570){\vector(1,-1){110}}
   \put(140,540){$\sigma$}
 \end{picture}
\end{center}
such that
\[ \vec{y}\circ g = f\circ \vec{x},\;\;\;  \vec{y}'\circ g' = f\circ \vec{x}',\;\;\;{\rm and}\;\;\; A(g')(a') = A(g)(a) \]
Then we have
\[ \vec{y}\circ g =f\circ \vec{x}\circ \sigma = \vec{y}\circ g'\circ \sigma \]
As $\vec{y}$ is mono $g=g'\circ \sigma$. Thus
\[ A(g')(a') = A(g)(a)= A(g')(A(\sigma)(a)) \]

As $\tau$ is a cartesian natural transformation we get, from the fact that the naturality square for g' is a pullback, that $a'=A(\sigma)(a)$. But this means that
\[ (\vec{x}',a')\sim (\vec{x}',A(\sigma)(a))\sim  (\vec{x}'\circ\sigma,a)\sim (\vec{x},a) \]
i.e. $[\vec{x}',a']=[\vec{x},a]$ and $\hat{\tau}$ is cartesian.

To show the converse assume that $\hat{\tau}$ is a cartesian natural transformation and $f:(n]\ra (m]$ is a surjection. We need to show
that the weak pullback
\begin{center} \xext=600 \yext=550
\begin{picture}(\xext,\yext)(\xoff,\yoff)
\setsqparms[1`1`1`1;600`400]
 \putsquare(0,50)[A(n]`B(n]`A(m]`{B(m]};\tau_n`A(f)`B(f)`\tau_m]
 \end{picture}
\end{center}
satisfies the uniqueness condition, as well. Fix $a,a'\in A(n)$ such that
\[ A(f)(a)=A(f)(a'),\;\;\;{\rm and}\;\;\; \tau_n(a)=\tau_n(a') \]
By assumption the square
\begin{center} \xext=600 \yext=550
\begin{picture}(\xext,\yext)(\xoff,\yoff)
\setsqparms[1`1`1`1;600`400]
 \putsquare(0,50)[\hat{A}(n]`\hat{B}(n]`{\hat{A}(m]}`{\hat{B}(m]};\hat{\tau}_{(n]}`\hat{A}(f)`\hat{B}(f)`\hat{\tau}_{(m]}]
 \end{picture}
\end{center}
is a pullback. We have elements $[1_{(n]},a], [1_{(n]},a']\in \hat{A}(n]$ such that
\[ \tau_{(n]}([1_{(n]},a])  =[1_{(n]},\tau_n(a)])=[1_{(n]},\tau_n(a')])=  \tau_{(n]}([1_{(n]},a'])\]
and
\[ \hat{A}(f)([1_{(n]},a])=[1_{(n]},A(f)(a)]=[1_{(n]},A(f)(a')]=\hat{A}([1_{(n]},a']) \]
Thus $[1_{(n]},a]=[1_{(n]},a']$ and $a=a'$, i.e. $\tau$ is a cartesian natural transformation. $\boxempty$

\noindent
Instytut Matematyki,
Uniwersytet Warszawski,
ul. S.Banacha 2,
00-913 Warszawa, Poland\\
szawiel@mimuw.edu.pl,
zawado@mimuw.edu.pl

\end{document}